\documentclass[11pt]{article}

\usepackage{amsmath,amssymb,amsthm,enumerate,framed,mdwlist,a4wide}
\usepackage{color,colortab}
\usepackage{pgf,tikz}
\usepackage{mathrsfs}
\usetikzlibrary{arrows}
\usepackage[]{todonotes} 
\usepackage{booktabs}
 \usepackage{listings} 
\usepackage{multicol}
\usepackage{color}
\usepackage{xspace}

\usepackage{booktabs}
\usepackage{multirow}
\numberwithin{table}{section}

\numberwithin{figure}{section}

\usepackage{amssymb}
\usepackage{latexsym}
\usepackage{amsfonts}
\usepackage{amsmath}
\usepackage{mathtools}

\usepackage{longtable}
\usepackage{rotating}
\usepackage{lscape}
\usepackage{algpseudocode}
\usepackage{algorithm}
\usepackage{url}
\usepackage[mathscr]{euscript}
\usepackage{hyperref}
\hypersetup{
	colorlinks,
	citecolor=blue,
	filecolor=blue,
	linkcolor=blue,
	urlcolor=black
}
\usepackage{pifont} 

\usepackage[]{todonotes} 

\newtheorem{proposition}{Proposition}
\newtheorem{assumption}{Assumption}

\newtheorem{theorem}{Theorem}

\DeclareMathOperator*{\argmax}{argmax} 

\usepackage{mathrsfs}

\newcommand{\Sy}{{\mathscr S}_n}
\newcommand{\Syl}{{\mathscr S}_{n+\ell}}
\newcommand{\sdp}{{\mathscr S}_n^{+}}
\newcommand{\sdpp}{{\mathscr S}_n^{++}}
\newcommand{\sdpl}{{\mathscr S}_{n+\ell}^{+}}
\newcommand{\sdn}{{\mathscr S}_n^{-}}

\newcommand{\dotp}[2]{\langle #1, #2 \rangle}
\newcommand{\Aat}{{\mathscr A}^{\top}}

\newcommand{\Aa}{{\mathscr A}}

\newcommand{\RR}{\mathbb{R}}

\newcommand{\ppp}[1]{\tilde{#1}}

\newcommand{\Diag}{{\rm Diag}}

\newcommand{\trace}{{\rm trace}}

\newcommand{\ignore}[1]{}
\newcommand\R{{\mathbb R}}

\def\ignore#1{}
 
\newcommand{\korange}{\textcolor[rgb]{0,0,0}}

\algnewcommand{\Input}[1]{%
    \textbf{Input: }#1}

\allowdisplaybreaks

 \setbox\strutbox=\hbox{\vrule height7pt depth2pt width0pt}

\def\square{{\setbox0=\hbox{X}\hbox to \ht0{\vrule\hss\vbox to \ht0{
  \hrule width \ht0\vfil\hrule width \ht0}\vrule}}}

\title{Dealing with inequality constraints in \korange{large-scale} semidefinite relaxations 
	for graph coloring and maximum clique problems}

\author{Federico Battista\thanks{Department of Industrial and Systems Engineering, Lehigh University, Bethlehem PA 18015, USA.\\ E-mail: {\tt feb223@lehigh.edu}}  \and Marianna De Santis\thanks{Dipartimento di Ingegneria Informatica Automatica e Gestionale, Sapienza Università di Roma, Via Ariosto, 25 00185, Roma, Italy. E-mail: {\tt mdesantis@diag.uniroma1.it}}}

\begin{document}


\maketitle

\begin{abstract}
Semidefinite programs (SDPs) can be solved in polynomial time by interior point
methods. However, when the dimension of the problem
gets large, interior point methods become impractical \korange{in terms of both} computational
time and memory requirements. \korange{Certain first-order} methods, such as Alternating Direction
Methods of Multipliers (ADMMs), \korange{established as} suitable algorithms to deal with
\korange{large-scale} SDPs and gained growing attention \korange{over} the past decade.
In this paper, we focus on an ADMM designed for SDPs in standard form and extend 
it to deal with inequalities when solving SDPs in general form. 
Beside numerical results on randomly generated instances, where we show that our method compares favorably with respect to the state-of-the-art solver~\texttt{SDPNAL+}~\cite{yang2015sdpnal}, we present results on instances from SDP relaxations
of classical combinatorial problems such as the graph coloring problem 
and the maximum clique problem.
\korange{Through extensive numerical experiments,} we show that even an inaccurate dual solution, obtained at a generic iteration of our \korange{proposed} ADMM, can represent an \korange{efficiently recovered valid bound} on the optimal solution of the combinatorial problems considered, as long as an appropriate post-processing procedure is applied.
\end{abstract}

\par\smallskip\noindent
{\bf Keywords:} Semidefinite programming, Graph coloring problem, Maximum clique problem
\par\smallskip\noindent
{\bf MSC:} 90C22, 90C27, 90C06

\section{Introduction}
Interest on semidefinite programming has considerably grown during the last two decades and \korange{is partly} due to the fact 
that many practical problems in operations research and combinatorial optimization can be modeled or 
approximated by semidefinite programs~\cite{LaRe:05}. \korange{The purpose of this paper is} to focus on the use of augmented Lagrangian methods for dealing with semidefinite programming relaxations of two well-known 
combinatorial problems: the graph coloring problem and the maximum clique problem.
Augmented Lagrangian methods are known to be an alternative 
to interior point methods and currently represent the most popular 
first-order algorithms used to handle large-scale semidefinite programs~\cite{BuMo:03,BuMo:05,PoReWi:06,MaPoReWi:09}.
\korange{As a variant of augmented Lagrangian methods, Alternating Direction Methods of Multipliers (ADMMs) have gained increasing attention in recent years~\cite{Sun2015AC3,yang2015sdpnal,Wen2010,deSaReWie:2018, battista2023semidefinite}. 
Falling under the class of first-order methods, their success can be attributed to the avoidance of computation, storage, and factorization of large Hessian matrices. This in turn enables a significant increase in scalability compared to interior point methods. On the other hand, this comes at some cost to accuracy, which should be properly addressed in the scenario where the semidefinite problem is a} 
relaxation of \korange{some} combinatorial problem \korange{with the goal of} obtaining a valid bound on its optimal solution.
In order to overcome this issue, \korange{safe-bounding procedures have been recently developed (see e.g.~\cite{battista2023semidefinite,JaChayKeil2007,cerulli2020improving,wiegele2022sdp}). These methods have mostly been employed to recover a posteriori the inaccuracies that are produced by solvers and provide a valid bound on the optimum.}

\korange{In this paper, we begin by extending both an existing ADMM and a safe-bounding procedure to deal with SDPs with both equality and inequality constraints. Subsequently, we compare their performance with SDPNAL+~\cite{yang2015sdpnal}, an established state-of-the-art solver for large-scale SDPs that was awarded the Beale-Orchard-Hays Prize in 2018.
In particular, we present numerical experiments on randomly generated SDPs and on instances of well-known SDP relaxations from the literature for the maximum clique problem and the graph coloring problem. The goal is twofold. First, we demonstrate the robustness of our extended ADMM equipped with a safe-bounding procedure compared to a state-of-the-art solver. Second, we show that even a low-precision dual solution obtained during a generic iteration of our ADMM can serve as a valid and efficiently recovered bound on the optimal solution of the combinatorial problems considered. As a byproduct, an extensive collection of noteworthy SDP bounds for two fundamental combinatorial optimization problems are presented.}

\subsection{Notation and outline}
Let $\Sy$ be the set of $n$-by-$n$ symmetric matrices. \korange{Further,} let $\sdp~\subset~\Sy$ \korange{($\sdpp~\subset~\Sy$)} be the set of positive semidefinite \korange{(positive definite)} matrices and $\Sy^{-}~\subset~\Sy$ be the set of negative semidefinite matrices.
In the following, we denote by $\left\langle X,Y\right\rangle = \trace(XY)$ the standard inner product in $\Sy$.
Whenever a norm is used, we consider the Frobenius norm in the case of matrices and the Euclidean norm 
in the case of vectors.
\korange{Letting $M\in \RR^{m\times n}$, we denote by $\text{vec}(M)$ the $mn$-dimensional vector formed by stacking the
columns of $M$ on top of each other ($\text{vec}^{-1}$ is the inverse operation). Letting $v \in \RR^n$, we denote by $\Diag{(v)}$ the diagonal matrix having the elements of
$v$ on the main diagonal.}
\korange{We denote by $e_i$} the $i$-th vector of the standard basis in $\RR^n$. 
\korange{Given $S \in \Sy$, we denote by $(S)_{+}$ and $(S)_{-}$ the projections of $S$ onto
the positive semidefinite and negative semidefinite cones, respectively.}
Moreover we denote by $\lambda(S)$ the vector of the eigenvalues of $S$ and by $\lambda_{\min}(S)$ and 
$\lambda_{\max}(S)$ the smallest and largest eigenvalue of $S$, respectively. 
\korange{At last, we denote by $\mathbf{0}_{n}$ and $\mathbf{0}_{n\times \ell}$ the all-zero column vector of size $n$ and the all-zero $n \times \ell$ matrix, respectively.}

The paper is organized as follows.
In Section~\ref{sec:admm}, \korange{we succinctly review the ADMM algorithm \texttt{ADAL} for solving SDPs in standard form, originally presented in~\cite{PoReWi:06,MaPoReWi:09,Wen2010}.
In Section~\ref{sec:ext_admm}, \korange{we discuss} SDPs with both equality and inequality constraints and we extend the aforementioned method to handle such problems. \korange{Starting from a low-precision dual solution, we then apply} the procedure outlined in~\cite{cerulli2020improving} to our context to recover a valid dual bound on the 
optimal primal value. In Section~\ref{sec:relax}, drawing from existing literature, we review well-known SDP formulations for relaxations of two fundamental combinatorial optimization problems: the maximum clique and the graph coloring problems. These formulations, in addition to randomly generated SDPs, constitute the test set for the numerical experiments we report in Section~\ref{sec:numres}. At last, some conclusions are drawn in Section~\ref{sec:conc}.}

\korange{\section{ADAL: an ADMM for SDPs in Standard Form} \label{sec:admm}} 
\korange{
In this section we review the basic concepts of \texttt{ADAL}~\cite{PoReWi:06,MaPoReWi:09,Wen2010}, 
an Alternating Direction Method of Multipliers (ADMM) to address SDPs in standard form:
\begin{align}
\label{psdp:eq}
\min \,       & \left\langle C,X\right\rangle \nonumber\\
\mbox{ s.t. } &  \langle A^j, X \rangle = b_j,	\quad \forall j=1,\dots, m \\
              & X\in \sdp, \nonumber
\end{align}
where $C\in\Sy$, $A^j\in\Sy$, for all $j=1,\ldots,m$, and $b := (b_1,\dots ,b_m)\in \mathbb{R}^m$. By defining the linear operator $\Aa: \Sy \rightarrow \mathbb{R}^m$, with
$(\Aa X)_j = \left\langle A^j,X\right\rangle$ for $A^j\in\Sy$, Problem~\eqref{psdp:eq} can be rewritten as
\begin{align}
\label{sdp:psdp:eq}
\min         & \quad \left\langle C, X \right\rangle \nonumber\\
\mbox{ s.t.} & \quad \Aa X = b \tag{PSDP-ST}\\
             & \quad X\in \sdp. \nonumber
\end{align}\quad
The dual of \eqref{sdp:psdp:eq} is defined as 
\begin{align}
\label{sdp:dsdp:eq}
\max          &\quad b^T y \nonumber\\
\mbox{ s.t. } &\quad  { {\mathscr A}}^{\top} y + Z = C \tag{DSDP-ST}\\
              &\quad  Z\in \sdp, \nonumber
\end{align}
where ${ {\mathscr A}}^{\top}: \RR^m \rightarrow \Sy$ is the adjoint operator of $ \Aa$, namely ${ {\mathscr A}}^{\top} y = \sum_{i=1}^m y_i  A^i$
for $y \in \RR^m$.}

\korange{\begin{assumption}
    Slater condition holds for~\eqref{sdp:psdp:eq} and~\eqref{sdp:dsdp:eq}; that is, there exist matrices $\tilde{X}, \tilde{Z} \in \sdpp$, and $\Tilde{y} \in \RR^{m}$ satisfying $\Aa \tilde{X} = b$ and $\Aa^{\top} \Tilde{y} + \tilde{Z} = C$.
\end{assumption}}

\korange{Under this assumption, it is well known that strong duality holds. \korange{Hence, 
the following KKT conditions are necessary and sufficient for the optimality of
a triplet $(X, y, Z)$ (see, e.g.,~\cite[Sec. 4.2]{NeNe:94})}:
\begin{align}
\label{eq:optCond}
   \Aa  X &= b, \quad& { {\mathscr A}}^{\top} y + Z &=  C,\quad&    Z X &= 0,  \quad&  X &\in \sdp, \quad&  Z &\in \sdp. \notag
\end{align}}
The method we consider is based on the maximization of the augmented Lagrangian built over the dual problem.
Let $X\in \Sy$ be the Lagrange multiplier for the dual equation 
${{\mathscr A}}^{\top} y + Z - C = 0$ and $\sigma>0$ be fixed. 
The augmented Lagrangian of~\eqref{sdp:dsdp:eq} is defined as
$$
L_\sigma(y, Z; X) = b^Ty - \langle {{\mathscr A}}^{\top} y + Z - C, X\rangle -
\frac{\sigma}{2}\|{ {\mathscr A}}^{\top} y + Z - C\|^2.$$
In augmented Lagrangian methods applied to the dual~\eqref{sdp:dsdp:eq}
the problem  
\begin{equation}\label{eq:maxLagZ}
 \begin{aligned}
      &\max  \, && L_\sigma(y,Z; X) \\
      &\mbox{ s.t. } && y\in \RR^m, \quad Z\in \sdp,
\end{aligned}
\end{equation}
is addressed at every iteration, where $X$ is fixed and $\sigma>0$ is a penalty parameter.
\korange{When the maximization of the augmented Lagrangian $L_\sigma(y,Z; X)$ is performed by iteratively optimizing with respect to $y$ at first, and then with respect to $Z$, we are considering the ADMM originally proposed in~\cite{PoReWi:06,MaPoReWi:09} and then extended in~\cite{Wen2010}.}

\korange{We present this method following~\cite{Wen2010}, in which it is referred as \texttt{ADAL} (Alternating Direction Augmented Lagrangian)}.
At each iteration of \texttt{ADAL}, the new point $(y^{k+1}, Z^{k+1},X^{k+1})$ is computed by the following steps:
\begin{align} \label{eq:opty}
 y^{k+1} &= \argmax_{y\in \RR^m} L_{\sigma^k}(y, Z^{k}; X^{k}), \\
 \label{eq:optZ}
 Z^{k+1} &= \argmax_{Z\in \sdp} L_{\sigma^k}(y^{k+1}, Z; X^{k}),\\
 \label{eq:updateX}
 X^{k+1} &=  X^{k} + \sigma^k ( \Aat y^{k+1} + Z^{k+1} - C).
\end{align}
The update of $y$ in~\eqref{eq:opty} 
can be performed in closed form, as it
derives from the first-order optimality conditions of the problem on the right-hand side of~\eqref{eq:opty}\korange{, i.e., }
$y^{k+1}$ 
is the unique solution of
\[
\nabla_y L_{\sigma^k} (y, Z^{k};  X^{k}) = b -  \Aa(  X^k + \sigma^k ({ {\mathscr A}}^{\top} y +  Z^k -  C)) = 0.
\]
That is,
\[
y^{k+1}= ( \Aa { {\mathscr A}}^{\top})^{-1}\Big(\frac{1}{\sigma^k} b -  \Aa(\frac{1}{\sigma^k}  X^k +  Z^k -  C)\Big).
\]
 The update of $ Z$ in~\eqref{eq:optZ} is conducted by considering the equivalent problem
\begin{equation}\label{eq:proj}
 \min_{ Z\in \sdp} \| Z + W^{k+1}\|^2,
\end{equation}
where 
\[W^{k+1} = \frac{ X^k}{\sigma^k} -  C + { {\mathscr A}}^{\top} y^{k+1} .\]
Solving problem \eqref{eq:proj}, is equivalent to \korange{projecting} $W^{k+1}\in \Sy$
onto the (closed convex) cone $\sdn$ and \korange{computing} its additive inverse (see Algorithm~\ref{alg:ADAL}). 
Such a projection is computed via the spectral decomposition of the matrix $W^{k+1}$.
Finally, it is \korange{clear} to see that the update of $ X$ in~\eqref{eq:updateX} can be \korange{obtained as follows:}
\begin{align*}
  X^{k+1} & =  X^k + \sigma^k ({ {\mathscr A}}^{\top} y^{k+1} +  Z^{k+1}  -  C ) = \\
 & = \sigma^k( X^k/\sigma^k -  C + { {\mathscr A}}^{\top} y^{k+1} - ( X^k/\sigma^k -  C + { {\mathscr A}}^{\top} y^{k+1})_-) = \\
 & = \sigma^k( X^k/\sigma^k -  C + { {\mathscr A}}^{\top} y^{k+1} )_+ \korange{= \sigma^k(W^{k+1})_+}.
\end{align*}
We report in Algorithm~\ref{alg:ADAL} the scheme of \texttt{ADAL}.
\begin{algorithm}[ht]
    \caption{Scheme of \texttt{ADAL} from~\cite{Wen2010}}
    \label{alg:ADAL}
    \begin{algorithmic}[1]
        \State Choose $\sigma>0$, $\varepsilon >0$, ${X}\in \sdp$, $ Z\in \sdp$
        \State $\delta = \max \{ r_P, r_D\}$
        \While{$\delta > \varepsilon$}
        \State $y= ( \Aa \Aat)^{-1}\Big(\frac{1}{\sigma} b -  \Aa(\frac 1 
        \sigma  X - C +  Z)\Big)$ \label{yupdate}
        \State $W = X/\sigma - C + { {\mathscr A}}^{\top} y$
        \vskip0.1truecm
        \State $ Z = -W_-$ and $  X 
        =\sigma W_+$ \label{proj}
        \State $\delta = \max \{ r_P, r_D\}$
        \State Update $\sigma$ 
        \EndWhile
    \end{algorithmic}
\end{algorithm}
The method stops as soon as the following errors related to primal feasibility ($ \Aa  X = b$, $ X\geq 0$) and dual feasibility (${ {\mathscr A}}^{\top} y +  Z +  S =  C$) are below a certain \korange{threshold defined, respectively, as the following:}
\begin{equation*}
    \begin{aligned}
 r_P&= \frac{\| \Aa  X - b\|}{1+\|b\|}, \quad \quad &
 r_D&= \frac{\|{ {\mathscr A}}^{\top} y +  Z -  C\|}{1+ \| C\|}.
    \end{aligned}
\end{equation*}
More precisely,
the algorithm stops as soon as the quantity $\delta = \max \{ r_P, r_D\}$
is less than a fixed precision $\varepsilon>0$.
\korange{It should be noted that} the optimality conditions $ X\in \sdp$, $ Z\in \sdp$ \korange{and} $ Z X=0$
are satisfied up to machine accuracy throughout the algorithm thanks to the projections employed in \texttt{ADAL}.
\korange{In the convergence analysis proposed in~\cite{Wen2010}, Algorithm~\ref{alg:ADAL} is interpreted as a fixed point method, i.e., at each iteration of~\texttt{ADAL}, the update of the primal and dual variables $(X, Z)$ is the result of the combination of two operators, both of which are proven to be non-expansive.
Indeed, we have that $(X^{k+1}, Z^{k+1}) = \mathcal{P}(W(X^k, Z^k))$, where
$\mathcal{P}$ denotes the projection performed at Step~\ref{proj} in Algorithm~\ref{alg:ADAL}
and $W(X^k,Z^k) = X^k/\sigma^k - C + \Aat y(Z^k,X^k)$, being the update performed at Step~\ref{yupdate} in Algorithm~\ref{alg:ADAL}, i.e.,
\[y(Z^k, X^k) = 
( \Aa \Aat)^{-1}\Big(\frac{1}{\sigma^k} b -  \Aa(\frac 1\sigma^k  X^k - C +  Z^k)\Big).\]}

\korange{Using the non-expansivity of $\mathcal{P}$ and $W$ (see Lemma~3 and Lemma~4 in~\cite{Wen2010}, respectively), it is possible to prove the following result (see Theorem 2 in~\cite{Wen2010}): 
\begin{theorem}
The sequence $\{({X}^k,y^k,{Z}^k)\}$ generated by Algorithm~\ref{alg:ADAL} applied to Problem~\eqref{sdp:psdp:eq} from any starting
point $({X}^0,y^0,{Z}^0)$ converges to a solution $({X}^*,y^*,{Z}^*)\in \Omega^*$, where $\Omega^*$ is the set of primal and dual solutions of~\eqref{sdp:psdp:eq} and~\eqref{sdp:dsdp:eq}.
\end{theorem}
}

\section{\korange{\texttt{ADAL-ineq}: applying \texttt{ADAL} to SDPs in general form}}\label{sec:ext_admm}
\korange{The aim of our work is to address SDPs that include linear inequality constraints. In order to do so, we can still use \texttt{ADAL}, namely Algorithm \ref{alg:ADAL}, by applying it to the reformulation with only equality constraints obtained through the introduction of slack variables. However, when dealing with large-scale SDPs, the memory required by \texttt{ADAL} to solve such a reformulation would grow substantially.
In this section, we show how to rewrite the steps of \texttt{ADAL} in terms of the matrices defining the original problem in a more efficient way that will reduce the memory requirements. For clarity, we will refer to this new version of \texttt{ADAL} as to \texttt{ADAL-ineq}}.
\korange{Recall that an} SDP in general form can be formulated as follows:
\begin{equation}
\label{psdp}
\begin{aligned}
&\min \,       && \left\langle C,X\right\rangle\\
&\mbox{ s.t. } &&  \langle A^i, X \rangle\leq b_i,	&&\quad \forall i=1,\dots,\ell	\\
&              &&  \langle A^j, X \rangle = b_j,	&&\quad \forall j=\ell + 1,\dots, m \\
&              && X\in \sdp,
\end{aligned}
\end{equation}
where $C\in\Sy$, $A^i\in\Sy$, \korange{for all $i=1,\ldots,m$,} and $b\in\RR^{m}$.
A standard way to deal with \korange{inequalities in} problem \eqref{psdp} is to add slack variables $s_i\geq 0$, \korange{for all $ i=1,\dots,\ell$,} and expand the matrix variable $X$ to $\bar X \in \Syl$:
$$
\bar X:= \begin{pmatrix}
X & \mathbf{0}_{n\times\ell}\\
\mathbf{0}_{\ell\times n} &  \Diag{(s)}\\
\end{pmatrix}.
$$
Recall that if $B$ is a diagonal matrix, the constraint $B \succeq 0$ boils down to $B \geq 0$.
In particular, imposing $\bar X \succeq 0$ 
is equivalent to \korange{considering} $X \succeq 0$ and  $s \geq 0$. 
\korange{By expanding the matrices $A^i, \;  A^j,$ and $C$, to $\bar A^i$, $\bar A^j$ and $\bar C$, respectively, for all $i=1,\dots,\ell$ and $j=\ell+1,\dots,m$, as}
$$
\bar A^i:= \begin{pmatrix}
	A^i & \mathbf{0}_{n\times\ell}\\
	\mathbf{0}_{\ell\times n} &  e_i^T e_i\\
\end{pmatrix}, \quad \quad
\bar A^j:= \begin{pmatrix}
A^j & \mathbf{0}_{n\times\ell}\\
\mathbf{0}_{\ell\times n} &  \mathbf{0}_{\ell\times \ell}\\
\end{pmatrix}, \quad \quad 
\bar C:= \begin{pmatrix}
C & \mathbf{0}_{n\times\ell}\\
\mathbf{0}_{\ell\times n} &  \mathbf{0}_{\ell\times \ell}\\
\end{pmatrix},
$$
problem \eqref{psdp} can be rewritten as an SDP in standard form as follows: 
\begin{equation}
\label{sdp:standard1}
\begin{aligned}
&\min \,       && \left\langle \bar C,\bar X\right\rangle\\
&\mbox{ s.t. } &&  \bar\Aa \bar X = b \\
&              && \bar X\in \sdpl,
\end{aligned}
\end{equation}
where $b := (b_1,\dots ,b_m)\in \mathbb{R}^m$ and $\bar \Aa: \Syl \rightarrow \mathbb{R}^m$
is the linear operator 
$(\bar \Aa X)_i = \left\langle \bar A^i,X\right\rangle$ with $\bar A^i\in\Syl$, \korange{for all} $i=1,\ldots,m$. 
The dual problem of \eqref{sdp:standard1} is defined as 
\begin{equation}
\label{sdp:standarddual1}
\begin{aligned}
&\min \,       && b^T y\\
&\mbox{ s.t. } &&  {\bar {\mathscr A}}^{\top} y + \bar Z = \bar C \\
&              && \bar Z\in \sdpl,
\end{aligned}
\end{equation}
where ${\bar {\mathscr A}}^{\top}: \RR^m \rightarrow \Syl$ is the adjoint operator of $\bar \Aa$, namely ${\bar {\mathscr A}}^{\top} y = \sum_{i=1}^m y_i \bar A^i$
for $y \in \RR^m$.
Note that the matrix $\bar Z\in \Syl$ is a ``surplus" matrix variable that can be written as
\[\bar Z:= \begin{pmatrix}
Z & \mathbf{0}_{n\times\ell} \\
\mathbf{0}_{\ell\times n} & \Diag{(p)} \\
\end{pmatrix},\] with $p \in \mathbb{R}^\ell$.
In particular, the equality constraint in \eqref{sdp:standarddual1}
can be rewritten as 
\korange{
\[
\bar C - \bar \Aa^T ( y) - \bar{Z} = \begin{pmatrix}
C - \Aa^Ty - Z & & \mathbf{0}_{n\times\ell}\\
\mathbf{0}_{\ell\times n} & & \Diag{(-y - p)}\\
\end{pmatrix} = 0.
\]}
As for the SDPs in standard form, if we assume that both the primal~\eqref{sdp:standard1} and the dual~\eqref{sdp:standarddual1} problems have strictly feasible 
points (i.e. Slater's condition is satisfied) strong duality holds and $(y, \bar Z, \bar X)$ is optimal 
for~\eqref{sdp:standard1} and~\eqref{sdp:standarddual1} if and only if the following KKT conditions hold \korange{(see, e.g.,~\cite[Sec. 4.2]{NeNe:94})}:
\begin{equation}\label{eq:optCond}
\begin{aligned}
  \bar \Aa \bar X &= b, \quad& {\bar {\mathscr A}}^{\top} y +\bar Z &= \bar C,\quad&   \bar Z\bar X &= 0,  \quad& \bar X &\in \sdpl, \quad& \bar Z &\in \sdpl. \\
\end{aligned}
\end{equation}
In the following, we assume that the constraints formed through the operator $\bar\Aa$ are linearly independent.

The memory required to store the augmented matrices $\bar C, \bar A_i, \bar A_j, \bar Z$ and $\bar X$ \korange{substantially increases} with the number $\ell$ of inequalities 
and even using efficient sparse matrix implementations may be insufficient to computationally deal with large-scale problems.
\korange{Thus, we propose rewriting} the steps of~\texttt{ADAL} \korange{applied to Problem \eqref{sdp:standard1}} in terms of the original matrices $C, A^i$, and $X$, 
so that one \korange{is only required to store} the matrices that are actually defining the problem.
\korange{
\begin{proposition}\label{prop:ypdate}
Step~\ref{yupdate} in Algorithm~\ref{alg:ADAL} applied to Problem~\eqref{sdp:standard1}, reading as 
\[
y= ( \bar\Aa \bar\Aa^\top)^{-1}\Big(\frac{1}{\sigma} \bar b -  \bar\Aa(\frac 1 
        \sigma  \bar X - \bar C +  \bar Z)\Big)
\]
can be performed in terms of the matrices of Problem~\eqref{psdp} only, namely in terms of the matrices
$C, A^i$, for all $i = 1,\ldots,m$, $X$ and $Z$.  
\end{proposition}
\begin{proof}
Let $1 \leq i \leq \ell$ be a generic index of an inequality constraint and let $\ell + 1 \leq j \leq m$ 
be a generic index of an equality constraint, then the following holds:
\begin{alignat*}{3}
&&\dotp{\bar A^i}{\bar X} &&= \dotp{A^i}{X} &+ s_i, \\
&&\dotp{\bar A^j}{\bar X} &&= \dotp{A^j}{X} &, \\
&&\dotp{\bar C}{\bar X}    &&= \dotp{C}{X}  &.
\end{alignat*}
The linear map applied to $\bar X$ becomes:
$$\bar \Aa(\bar X) =
\begin{pmatrix}
\dotp{A^1}{X}\\
\vdots \\
\dotp{A^\ell}{X}\\
\dotp{A^{\ell+1}}{X}  \\
\vdots \\
\dotp{A^m}{X} \\
\end{pmatrix} + 
\begin{pmatrix}
s^T \\
\mathbf{0}_{m - \ell}\\
\end{pmatrix} = \Aa(X) +
\begin{pmatrix}
s^T \\
\mathbf{0}_{m - \ell}\\
\end{pmatrix}. 
$$
Similarly, the adjoint operator $\bar \Aa^T : \mathbb{R}^m \rightarrow \Syl$ of $\bar \Aa$ is defined as 
$$
\bar \Aa^Ty := \sum_{i=1}^m y_i\bar A^i = 
\begin{pmatrix}
\sum_{i=1}^m y_iA^i & &  \mathbf{0}_{n\times\ell}\\
   \mathbf{0}_{\ell\times n} & & \Diag{(y)}\\
\end{pmatrix} = 
\begin{pmatrix}
\Aa^Ty & & \mathbf{0}_{n\times\ell}\\
\mathbf{0}_{\ell\times n} & & \Diag{(y)}
\end{pmatrix}.
$$
Using the operator $\text{vec}{}$, we can write $\bar \Aa( \bar X) = b$ as $\bar A\ \text{vec}(\bar X) = b$, where
$$
\bar A := \left( \text{vec}(\bar A^1), \dots, \text{vec}( \bar A^m) \right) ^T \in \mathbb{R}^{m\times (n+\ell)^2}.
$$
Note that the matrix $\bar A^i$, for all $i=1,\dots,\ell$, corresponding to the $i$-th inequality constraint, is the unique matrix having $1$ in position $(n + i, n + i)$. 
Then, $\bar A \bar A^T$ can be expressed in terms of $AA^T$ as follows:
$$
\bar A \bar A^T = AA^T + 
\Diag{
\begin{pmatrix}
\mathbf{1}_\ell \\
\mathbf{0}_{m-\ell}
\end{pmatrix}},
$$
where $A := \left( \text{vec}(A^1), \dots, \text{vec}(A^m) \right)^T \in \mathbb{R}^{m\times n^2}.$
Indeed, the zero entries of $\bar A^i$ do not contribute in the row-by-column product and
the $1$ in position $(n + i, n + i)$ contributes only to the entry where vec$(\bar A^i)$ is multiplied by itself, 
 i.e., in position $(i, i)$ of $\bar A \bar A^T$.
According to the notation introduced, the update of the $y$ variable can be rewritten as follows:
\[y^{k+1} = \left( AA^T + 
\Diag
\begin{pmatrix}
\mathbf{1}_\ell \\
\mathbf{0}_{m-\ell}
\end{pmatrix}\right)^{-1}
\left(\frac{1}{\sigma^k} b - A\ \text{vec}\left(\frac{1}{\sigma^k} X^k -  C + Z^k \right) + \begin{pmatrix}
\frac{1}{\sigma^k}{s^{k}}^T + {p^{k}}^T \\
\mathbf{0}_{m - \ell}\\
\end{pmatrix}  \right). 
\]
\end{proof}
}
%
\korange{As a consequence of Proposition~\ref{prop:ypdate}, the spectral decomposition of the matrix $W$, performed in Step~\ref{proj} and needed for updating the variables $\bar X$ and $\bar Z$, can also be computed without storing the augmented matrices.
Indeed $W^{k+1}$ can be written in a ``block-wise'' fashion as follows:
\[
W^{k+1}
= 
\begin{pmatrix}
	\frac{1}{\sigma^k} X^k - C + \Aat y^{k+1} & \mathbf{0}_{n\times\ell}\\
	\mathbf{0}_{\ell\times n}& \Diag\left( \frac{1}{\sigma^k}{s^k} + y^{k+1}\right)
\end{pmatrix}.
\]}

\korange{To compute the eigenvalues and eigenvectors of $W^{k+1}$, we begin by performing the spectral decomposition of the matrix $\frac{X^k}{\sigma^k} - C + \Aa^T y^{k+1}$. Then, we straightforwardly obtain the eigenvalues and eigenvectors corresponding to the diagonal part of $W^{k+1}$. Finally, we adjust the dimension of the computed eigenvectors to ensure they belong to $\R^{n+\ell}$.}
  
\korange{Proposition~\ref{prop:ypdate} and the reasoning reported above demonstrate that it is possible to apply 
Algorithm~\ref{alg:ADAL} to Problem~\eqref{sdp:standard1} without storing the augmented matrices
by rewriting its steps in terms of the matrices
that are defining the original inequality constrained problem~\eqref{psdp}.
As already mentioned, we denote this version of \texttt{ADAL} as \texttt{ADAL-ineq}
and its scheme is reported in Algorithm~\ref{alg:ADALin}.
\begin{algorithm}[ht]
    \caption{\korange{Scheme of \texttt{ADAL-ineq}}}
    \label{alg:ADALin}
    \begin{algorithmic}[1]
        \State Given Problem \eqref{sdp:standard1}, choose $\sigma>0$, $\varepsilon >0$, 
        ${X}\in \sdp$, $ Z\in \sdp$, $s\in \R^\ell$, $p\in \R^\ell$
        \State $\delta = \max \{ r_P, r_D\}$
        \While{$\delta > \varepsilon$}
        \State $y=  
         \left( AA^T + \Diag
            \begin{pmatrix}
                \mathbf{1}_\ell \\
                \mathbf{0}_{m-\ell}
            \end{pmatrix}\right)^{-1}
            \left(\frac{1}{\sigma} b - A\ \text{vec}\left(\frac{1}{\sigma} X -  C + Z \right) + \begin{pmatrix}
                \frac{1}{\sigma}s^T + p^T \\
                \mathbf{0}_{m - \ell}\\
            \end{pmatrix}  \right)
        $ \label{yupdatein}
        \vskip0.2truecm
        \State $W = 
        \begin{pmatrix}
	           \frac{1}{\sigma} X - C + \Aat y & \mathbf{0}_{n\times\ell}\\
	           \mathbf{0}_{\ell\times n}& \Diag\left( \frac{1}{\sigma}{s} + y\right)
        \end{pmatrix}
        $
        \vskip0.2truecm
        \State $ Z =-W_-$ and $X 
        =\sigma W_+$ \label{projin}
        \State $\delta = \max \{ r_P, r_D\}$
        \State Update $\sigma$ 
        \EndWhile
  \end{algorithmic}
\end{algorithm}
}

\korange{
From the convergence of \texttt{ADAL}~\cite{Wen2010} we can state the following:
\begin{theorem}
The sequence $\{(\bar{X}^k,y^k,\bar{Z}^k)\}$ generated by \texttt{ADAL-ineq} from any starting
point $(\bar{X}^0,y^0,\bar{Z}^0)$ converges to a solution $(\bar{X}^*,y^*,\bar{Z}^*)\in \Omega^*$, where $\Omega^*$ is the set of primal and dual solutions of~\eqref{sdp:standard1} and~\eqref{sdp:standarddual1}.
\end{theorem}
}
\subsection{Obtaining dual bounds}\label{sec:dualbound}
\korange{The approximation of combinatorial problems is one of the most relevant applications of semidefinite programming. This is because} the optimal solution 
of a semidefinite relaxation can be computed in polynomial time and generally gives a better bound than that obtained solving a linear relaxation (see e.g.~\cite{matRel}).
Given a pair of primal-dual SDPs, weak and strong duality hold under the assumption that both problems are strictly feasible. 
Duality results imply that the objective function value of every feasible
solution of the dual SDP is a valid bound on the optimal
objective function value of the primal. 
Therefore, every dual feasible solution, and in particular the optimal dual
solution of an SDP relaxation, gives a valid bound on the solution of the related combinatorial
optimization problem.
 Hence, being able to compute dual feasible solutions - even of moderate quality - can be extremely useful 
 when considering branch-and-bound frameworks to define exact solution methods
for specific combinatorial optimization problems.
 Following ideas developed in~\cite{cerulli2020improving}, we define a post-processing procedure 
for \texttt{ADAL-ineq} on general SDPs, that allows \korange{one to obtain} a feasible dual solution starting from a positive semidefinite matrix $\ppp{Z}\in \sdp$.
Let $\Aa_{ineq}$ and $\Aa_{eq}$ be the linear operators defining the inequality and equality constraints in 
problem~\eqref{psdp}\korange{, respectively; thus $\Aa_{ineq} = \left\langle A^i,X\right\rangle$ with $A^i\in\Sy$ and 
$\Aa_{eq} = \left\langle A^j,X\right\rangle$ with $A^j\in\Sy$, for all $i=1,\ldots,\ell$ and $j=\ell+1,\ldots,m$.}
Let $b_{ineq}$ and $b_{eq}$ be the \korange{right-hand-side} vectors accordingly defined.
Introducing the adjoint operators of $\Aa_{ineq}$ and $\Aa_{eq}$, the dual problem~\eqref{psdp} can be equivalently written as
\begin{equation}
\label{eq:dualNJ}
\begin{aligned}
&\max \,       && -b_{ineq}^T \lambda + b_{eq}^T\mu\\
&\mbox{ s.t. } &&  C + \Aat_{ineq} \lambda - \Aat_{eq} \mu  = Z \\
&              && Z\in \sdp,\; \lambda \geq 0,
\end{aligned}
\end{equation}
with $\lambda \in \R^\ell$ and $\mu \in \R^{m-\ell}$.
We can then extend the results proposed in~\cite{cerulli2020improving} and define a procedure to 
\korange{obtain} feasible solutions of problem~\eqref{eq:dualNJ} \korange{along with}, by weak duality, 
valid bounds on the optimal
objective function value of the primal~\eqref{psdp}.
Let $\ppp{Z} \in \sdp$.  If the linear programming problem 
\begin{equation}\label{eq:LP}
\begin{aligned}
&\max \,       && -b_{ineq}^T \lambda + b_{eq}^T\mu\\
&\mbox{ s.t. } &&  C + \Aat_{ineq} \lambda - \Aat_{eq} \mu  = \ppp{Z} \\
&              && \lambda \geq 0
\end{aligned}
\end{equation}
has an optimal solution $(\ppp{\lambda},\ppp{\mu}) \in \R^{m}$, 
then $(\ppp{\lambda}, \ppp{\mu}, \ppp{Z})$ is a feasible solution for~\eqref{eq:dualNJ}
and the value $-b_{ineq}^T \ppp{\lambda} + b_{eq}^T\ppp{\mu}$ \korange{yields} a dual bound.
If~\eqref{eq:LP} is unbounded, \korange{then is~\eqref{eq:dualNJ} also unbounded,
implying that} the primal~\eqref{psdp} is not feasible.
If~\eqref{eq:LP} is infeasible, \korange{then the initial $\ppp{Z} \in \sdp$ permits neither a feasible dual solution nor a dual bound.}
From a practical viewpoint, once problem~\eqref{psdp} is approximately solved by~\texttt{ADAL-ineq},
\korange{one can try to obtain} a feasible solution \korange{to} problem~\eqref{eq:dualNJ} by addressing problem~\eqref{eq:LP}.

\section{Bounding the clique number and the chromatic number of a graph}\label{sec:relax}
Given an undirected graph $G = (V,E)$, where $V$ is the set of vertices and $E$ is the set of edges,
\korange{a set $W\subseteq V$ is a \emph{clique} if every two vertices in $W$ are adjacent, while it is called \emph{stable} if no two vertices in $W$ are adjacent.} 
The \emph{clique number} $\omega(G)$ and the \emph{stability number} $\alpha(G)$ \korange{are the maximum cardinalities of a clique and a stable set in $G$, respectively.} 
A $k$-coloring is a partition of $V$ into $k$ stable sets. The \emph{chromatic number} $\chi(G)$ is the smallest integer $k$ for which $G$ has a $k$-coloring.
Denoting with $\bar G = (V,\bar E)$ the complementary graph of $G$, it holds \korange{that}
\[  \omega(\bar G) = \alpha(G)  \leq \chi(\bar G).\]
\korange{Lov\'asz}~\cite{lovasz1979shannon} introduced the so called \emph{theta number} $\vartheta (G)$ that is an upper bound \korange{on} the clique number $\omega(\bar G)$, the stability 
number $\alpha(G)$, and is \korange{also} a lower bound for the chromatic number $\chi(\bar G)$.
The important property of $\vartheta (G)$ is that it can be computed with an arbitrary precision in polynomial time, as it is 
the optimal value of the following SDP~\cite{grotschel2012geometric}:
\begin{equation}\label{eq:lovasztheta}
\begin{array}{l l l l}
\vartheta(G)  = & \max  & \left\langle J, X \right\rangle \\[1.1ex]
 &\mbox{ s.t. } & \trace(X) = 1 \\[1.1ex]
&  &  X_{ij} = 0 & \{i,j\} \in E\\[1.1ex]
  & & X\in \sdp, \nonumber
\end{array}
\end{equation}
where $J$ is the $n$-by-$n$ matrix of all ones, \korange{with} $n = |V|$. 
Starting from this relaxation, several attempts for sharpening $\vartheta (G)$ as a bound for $\omega(\bar G)$, $\alpha(G)$ and $\chi(\bar G)$ have been made (see, e.g.,~\cite{DUKANOVIC2008180,gruber2003computational,giandomenico2009application,giandomenico2013strong,giandomenico2015ellipsoidal,locatelli2015improving,GaarRendl}). 
As a first way to improve $\vartheta(G)$, we consider the numbers $\vartheta_+(G)$ and $\bar \vartheta_+(G)$
obtained as solutions of the following SDPs, where bounds on the entries of the matrix variables are introduced:
\begin{equation}
\begin{array}{l l l }
\vartheta_+(G)  = & \max  & \left\langle J, X \right\rangle \\[1.1ex]
 &\mbox{ s.t. } & \trace(X) = 1 \\[1.1ex]
&  &  X_{ij} = 0 \quad \{i,j\} \in \korange{E}\\[1.1ex]
  & & X\geq 0 \nonumber\\[1.1ex]
  & & X\in \sdp, \nonumber\\
\end{array}
\quad \quad \quad \quad 
\begin{array}{l l l l}
\bar\vartheta_+(G)  = & \max  & \left\langle J, X \right\rangle \\[1.1ex]
 &\mbox{ s.t. } & \trace(X) = 1 \\[1.1ex]
&  &  X_{ij} \leq 0 \quad \{i,j\} \in \korange{E}\\[1.1ex]
  & & X\in \sdp. \nonumber\\
\end{array}
\end{equation}
The values $\vartheta_+(G)$ and $\bar \vartheta_+(G)$ are related to 
$\omega(\bar G),\, \alpha(G)$ and $\chi(\bar G)$
as follows 
\[\omega(\bar G) = \alpha(G) \leq \vartheta_+(G) \leq \vartheta(G) \leq \bar\vartheta_+(G) \leq \chi(\bar G).\]
In the literature, equivalent formulations for both $\vartheta_+(G)$ and $\bar\vartheta_+(G)$ have been proposed~\cite{LaRe:05}
and for our computational \korange{experiments}, we consider the following formulation for $\bar\vartheta_+(G)$:
\begin{equation}
    \begin{array}{l l l l}
\bar\vartheta_+(G)  = & \min  & t \\[1.1ex]
 &\mbox{ s.t. } & X_{ii} = t-1 &  i\in \korange{V}\\[1.1ex]
 & & X_{ij} = -1 &  \{i,j\} \in \korange{\bar E}\\[1.1ex]
 & & X_{ij} \geq -1 &  \{i,j\} \in \korange{E}\\[1.1ex]
 & & X\in \sdp,\ t \in \R_+. \label{eq:theta_gc}
\end{array}
\end{equation}
\korange{where $t$ is an additional auxiliary variable. 
Problem~\eqref{eq:theta_gc} can be reformulated as~\eqref{psdp}
by including $t$ in the matrix variable as an additional
element on the diagonal.}
Note that, in both the formulations of $\vartheta_+(G)$ and $\bar\vartheta_+(G)$, 
the entries of the matrix $X$ are bounded from below.
In the context of ADMMs defined over the dual problem, 
bounds on the matrix variable can be handled by introducing 
a further step, where a projection onto the nonnegative orthant is
performed (see e.g.~\cite{Wen2010, cerulli2020improving, wiegele2022sdp}).
Although these 3-block ADMMs may not theoretically converge~\cite{chen2016direct}, they perform well in practice.

\section{Numerical results}\label{sec:numres}
\korange{In this section, we present the results of our computational study, where we evaluate the performance of both \texttt{ADAL-ineq} and \texttt{SDPNAL+}~\cite{yang2015sdpnal}. The comparison is conducted on randomly generated instance, as well as on SDP relaxations of both the stable set and the graph coloring problems.
The software \texttt{SDPNAL+}, accessible at \url{https://blog.nus.edu.sg/mattohkc/softwares/sdpnalplus/}, integrates an ADMM with a semismooth Newton-Conjugate Gradient method. It is implemented in MATLAB, using a refined management of the matrices exploiting their symmetry. This approach allows the optimization of a significant portion of its C subroutines, which are provided through Mex files.}
\korange{\texttt{ADAL-ineq} is developed in MATLAB utilizing its built-in functions, and is available at~\url{https://github.com/batt95/ADAL-ineq}, along with the instances used in our numerical experiments and an alternative Python implementation.}

\korange{The numerical performance of ADMMs, including \texttt{ADAL-ineq}, strongly depend on the update rule used for the penalty parameter $\sigma$. As in~\cite{cerulli2020improving,wiegele2022sdp},
we follow the strategy by Lorenz and Tran-Dinh~\cite{Lorenz2019}, considering at every iteration $k$ the ratio
between the norm of the primal variable $ X^k$ and norm of the dual variable $ Z^k$.
In the implementation of our post-processing procedure, described in Section \ref{sec:dualbound}, we 
used Gurobi 9.1.1~\cite{gurobi} as the solver for problem~\eqref{eq:LP}.}
The experiments were carried out on an Intel(R) Xeon(R) CPU E5-2698 v4 running at 2.20GHz, 
with 256GB of RAM, under Linux (Ubuntu 16.04.7).

We compare the performance of the algorithms using performance \korange{profiles 
proposed} by Dolan and Mor\'e~\cite{DM2002}.
Given a set of solvers~$\mathcal{S}$ and a set of problems $\mathcal{P}$, 
the performance of a solver $s \in \mathcal{S}$ on problem $p \in \mathcal{P}$ is compared
against the best performance obtained by any solver in~$\mathcal{S}$ 
on the same problem. The performance ratio is defined as
$
r_{p,s} = t_{p,s}/\min\{t_{p,s^\prime} \mid s^\prime \in\mathcal{S}\},
$
where $t_{p,s}$ is the measure we want to compare, and we consider a cumulative distribution
function~$\rho_s(\tau) = |\{p\in \mathcal{P} \mid r_{p,s}\leq \tau \}| /|\mathcal{P}|$.
The performance profile for $s \in S$ is the plot of the function $\rho_s$.

\subsection{Comparison on randomly generated instances}
The random \korange{SDPs} considered in the first experiment \korange{are created from an adaptation of} the
instance generator used in~\cite{MaPoReWi:09}. 
\korange{Given a triplet $(n, m, p) \in \mathbb{N} \times \mathbb{N} \times [0, 1]$, the resulting SDP consists of a matrix variable $X \in \sdp$ and includes $m$ linear constraints, of which $round(pm)$ are inequalities.
Instances are generated with values from $n \in \{200, 250, 500, 1000\}$, $m \in \{5000, 10000, 25000, 50000, 100000\}$ and $p \in \{0.25, 0.5, 0.75\}$. For each parameter combination, we generate 5 different instances, excluding values of $n$ and $m$ that result in a constraint matrix $A \in \RR^{m \times n^2}$ with linearly dependent rows. The final test set counts 150 random SDPs. A time limit of 1800 seconds of CPU time is set.}

\begin{figure}
    \centering
    \includegraphics[width=0.7\textwidth]{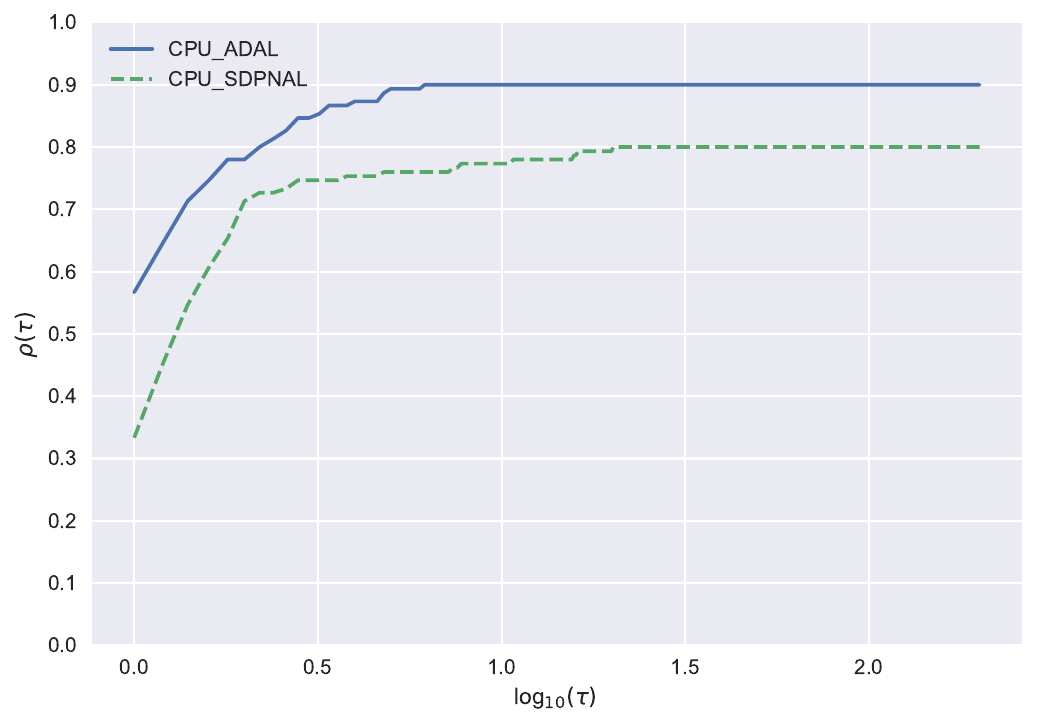}
    \caption{Performance profiles on CPU time. Comparison between \texttt{ADAL-ineq} and \texttt{SDPNAL+} on random instances.}
    \label{fig:rand}
\end{figure}

\korange{In Table~\ref{tab:rand}, we report the comparison between
\texttt{ADAL-ineq} and \texttt{SDPNAL+} in terms of number of iterations and CPU time needed in order to reach an accuracy of $10^{-5}$. For each solver and each combination of $n$, $m$ and $p$, we 
report the number of instances solved within the time limit 
along with the average running time.}
We notice that for $n = 250$ and $m = 25000$, \texttt{SDPNAL+} is not able to solve any instance within the time limit,
while \texttt{ADAL-ineq} is able to solve all of them with a precision of $10^{-5}$.
For $n = 500$ and $m = 100000$, both algorithms are not able to solve any instance within the time limit.
\texttt{SDPNAL+} performs better on instances with $n = 1000$ and $m = 10000$, while for the other instances either the two solvers show similar performances or \texttt{ADAL-ineq} outperforms \texttt{SDPNAL+}.
\begin{table}
	\centering
	{\scalebox{.85}{
\begin{tabular}{lll|cc|cc}
\toprule
  & & & \multicolumn{2}{c}{\texttt{ADAL-ineq}} & \multicolumn{2}{|c}{\texttt{SDPNAL+}}\\
\korange{$n$} & \korange{$m$} & \korange{$p (\%)$} & \#sol  & CPU time &  \#sol  &  CPU time  \\

\midrule
200  & 10000  & 25 & 5 &   39.24 & 5 &    33.05 \\
     &        & 50 & 5 &   58.24 & 5 &   109.14 \\
     &        & 75 & 5 &   67.14 & 5 &   713.82 \\
\toprule
250  & 5000   & 25 & 5 &    7.99 & 5 &    11.05 \\
     &        & 50 & 5 &    9.87 & 5 &    15.51 \\
     &        & 75 & 5 &   11.28 & 5 &    16.93 \\
\midrule
     & 25000  & 25 & 5 &  838.04 & 0 &  - \\
     &        & 50 & 5 & 1166.45 & 0 &  - \\
     &        & 75 & 5 & 1114.52 & 0 &  - \\
\toprule
500  & 10000  & 25 & 5 &   15.52 & 5 &    15.54 \\
     &        & 50 & 5 &   16.49 & 5 &    22.45 \\
     &        & 75 & 5 &   28.87 & 5 &    23.94 \\
\midrule
     & 25000  & 25 & 5 &   18.11 & 5 &    31.33 \\
     &        & 50 & 5 &   30.20 & 5 &    50.78 \\
     &        & 75 & 5 &   45.53 & 5 &    52.57 \\
\midrule
     & 50000  & 25 & 5 &  217.61 & 5 &   106.28 \\
     &        & 50 & 5 &  260.43 & 5 &   221.66 \\
     &        & 75 & 5 &  325.71 & 5 &   250.97 \\
\midrule
     & 100000 & 25 & 0 & 	- 	& 0 &  - \\
     &        & 50 & 0 & 	- 	& 0 &  - \\
     &        & 75 & 0 & 	- 	& 0 &  - \\
\toprule
1000 & 10000  & 25 & 5 &  136.63 & 5 &    49.52 \\
     &        & 50 & 5 &  157.21 & 5 &    58.22 \\
     &        & 75 & 5 &  242.63 & 5 &    71.38 \\
\midrule
     & 50000  & 25 & 5 &   57.19 & 5 &    60.96 \\
     &        & 50 & 5 &   94.09 & 5 &   109.48 \\
     &        & 75 & 5 &  110.00 & 5 &   111.29 \\
\midrule
     & 100000 & 25 & 5 &   83.15 & 5 &   136.53 \\
     &        & 50 & 5 &  127.37 & 5 &   181.13 \\
     &        & 75 & 5 &  155.05 & 5 &   184.21 \\
\bottomrule
\end{tabular}
}}
\label{tab:rand}
\caption{Results on \korange{$150$} random instances}
\end{table}
The performance profiles of \texttt{ADAL-ineq} and \texttt{SDPNAL+} on random instances are reported in Figure~\ref{fig:rand}, \korange{highlighting the superior} performance of \texttt{ADAL-ineq} with respect to \texttt{SDPNAL+}\korange{; on close to} $60\%$ of the instances \texttt{ADAL-ineq} is the fastest algorithm and is also able to solve $90\%$ of the instances \korange{whereas} \texttt{SDPNAL+} \korange{is only able to solve} $80\%$ of the instances within the time limit.

\subsection{Comparison on instances from SDP relaxations of the maximum clique problem}
\label{sec:res_clique}
\korange{We now} report the results on the SDP relaxation $\vartheta_+(G)$ for bounding the clique number (or the stability number) of a graph. 
We considered graphs from the second DIMACS implementation challenge~\cite{Johnson:1996}, available at \url{ftp://dimacs.rutgers.edu/pub/challenge/graph/benchmarks/clique}. These graphs form the standard benchmark for the maximum clique problem. \korange{Hence, we considered the complement} of these graphs to convert the maximum clique instances into stable set problem instances. \korange{For this experiment, a time limit of $3600$ seconds was set.}

\korange{In addition to employing \texttt{ADAL-ineq} and \texttt{SDPNAL+} to determine $\vartheta_+(G)$ 
and stopping the algorithms when the termination criteria were satisfied with a precision of $10^{-6}$, we equipped \texttt{ADAL-ineq} with the post-processing procedure \korange{detailed in Section \ref{sec:dualbound}}, applied every $200$ iterations and after termination.}
Every time the post-processing procedure is called, we give as input the matrix $Z^k$ obtained by~\texttt{ADAL-ineq} at the corresponding iteration $k$
and solve the linear programming problem~\eqref{eq:dualNJ}
using Gurobi~\cite{gurobi}.
Along \korange{with} the iterations of~\texttt{ADAL-ineq}, we \korange{also store} in memory 
the best dual bound found by the post-processing procedure,
together with the CPU time needed to detect it.
Note that every dual bound computed by the post-processing procedure and, in particular, the best dual bound \korange{are valid upper bounds} on the stability number.

\korange{Table~\ref{tab:theta_ss} is organized as follows: for each instance we report \korange{the} name (Graph), values of $\vartheta_+(G)$ arising from, respectively, the dual objective function values returned by \texttt{ADAL-ineq}, \texttt{SDPNAL+} and the best dual bound found (BestBound) throughout the iterations of~\texttt{ADAL-ineq}, along with the CPU times needed by~\texttt{ADAL-ineq} and~\texttt{SDPNAL+} to 
reach the stopping criterion, to identify the BestBound and the total time spent by the post-processing procedure.}
\korange{Note that the BestBound reported in the table is, in general, obtained in one of the post-processing calls when running \texttt{ADAL-ineq} and is not necessarily the one arising from the last call of the post-processing. We highlight in bold the values of the bounds found on the instances in which they differ.}

\korange{It should be noted} that the post-processing procedure applied \korange{within} \texttt{ADAL-ineq} is able to compute valid dual bounds on every instance but on \texttt{keller6}, where \texttt{ADAL-ineq} shows a failure. On \texttt{p\_hat1500-2}, even if both~\texttt{ADAL-ineq} and~\texttt{SDPNAL+} did not converge within $3600$ seconds, the post-processing procedure \korange{was} able to compute a valid dual bound. \korange{In order to understand the quality of this bound, we ran 
\texttt{SDPNAL+} on the \texttt{p\_hat1500-2} instance until convergence. Interestingly, it turns out that the BestBound and the optimal value match. By providing a valid dual bound, our procedure overcomes the impossibility of the solvers to reach termination criteria although the solution they achieve is close to an optimal one.}
\korange{Consider} that for huge graphs the time needed to find the best dual bound may be greater than the time needed by \texttt{ADAL-ineq} to converge\korange{. T}his comes from the fact that the best dual bound can be recovered at the last iteration computed.
\korange{We also wish to highlight} that the overall time needed to apply the post-processing procedure in~\texttt{ADAL-ineq} is small with respect to the overall time needed by the algorithm, and clearly it may be lowered by seldom applying the procedure.

\begin{table}
	\centering
	{\scalebox{.75}{
\begin{tabular}{l|rrr|rrrr}
\toprule
 &   \multicolumn{3}{c}{$\vartheta_+(G)$} &  \multicolumn{4}{c}{CPU times}  \\
 \hline
Graph           &     \texttt{ADAL-ineq}       &      \texttt{SDPNAL+}       &  BestBound         &                    \texttt{ADAL-ineq}       &      \texttt{SDPNAL+}       &  BestBound             &    post-proc     \\
\midrule
\texttt{DSJC125.1}   &      38.04 &       38.04 &     38.04 &          2.89 &            3.52 &           2.56 &    0.12 \\
\texttt{DSJC125.5}   &      11.40 &       11.40 &     11.40 &          1.93 &            0.77 &           1.70 &    0.09 \\
\texttt{DSJC125.9}   &       4.00 &        4.00 &      4.00 &          2.41 &            1.17 &           2.44 &    0.09 \\
\texttt{DSJC500-5}   &      22.57 &       22.57 &     22.57 &          6.83 &            4.85 &           7.30 &    0.48 \\
\texttt{DSJC1000-5}  &      31.67 &       31.67 &     31.67 &         41.60 &           34.32 &          43.69 &    3.59 \\
\texttt{C125-9}      &      37.55 &       37.55 &     37.55 &          2.88 &            0.99 &           2.73 &    0.11 \\
\texttt{C250-9}      &      55.82 &       55.82 &     55.82 &          7.82 &            3.12 &           7.15 &    0.43 \\
\texttt{C500-9}      &      83.58 &       83.58 &     83.58 &         30.55 &            8.32 &          31.03 &    1.48 \\
\texttt{C1000-9}     &     122.60 &      122.60 &    122.60 &        159.79 &           34.93 &         147.00 &    9.74 \\
\texttt{C2000-5}     &      44.56 &       44.56 &     44.56 &        389.59 &          534.67 &         398.86 &   22.42 \\
\texttt{C2000-9}     &     177.73 &      177.73 &    177.73 &       1238.53 &          278.60 &        1247.62 &   68.30 \\
\texttt{brock200\_1}  &      27.20 &       27.20 &     27.20 &          3.45 &            1.06 &           3.12 &    0.24 \\
\texttt{brock200\_2}  &      14.13 &       14.13 &     14.13 &          1.91 &            1.08 &           1.98 &    0.15 \\
\texttt{brock200\_3}  &      18.67 &       18.67 &     18.67 &          2.24 &            1.07 &           2.31 &    0.15 \\
\texttt{brock200\_4}  &      21.12 &       21.12 &     21.12 &          2.83 &            0.96 &           2.92 &    0.18 \\
\texttt{brock400\_1}  &      39.33 &       39.33 &     39.33 &          7.81 &            3.84 &           8.09 &    0.53 \\
\texttt{brock400\_2}  &      39.20 &       39.20 &     39.20 &          8.30 &            4.02 &           8.58 &    0.50 \\
\texttt{brock400\_3}  &      39.16 &       39.16 &     39.16 &          8.64 &            3.59 &           8.92 &    0.48 \\
\texttt{brock400\_4}  &      39.23 &       39.23 &     39.23 &          7.99 &            3.53 &           8.27 &    0.49 \\
\texttt{brock800\_1}  &      41.87 &       41.87 &     41.87 &         24.31 &           11.67 &          20.96 &    2.66 \\
\texttt{brock800\_2}  &      42.10 &       42.10 &     42.10 &         23.94 &           12.88 &          20.49 &    2.59 \\
\texttt{brock800\_3}  &      41.88 &       41.88 &     41.88 &         24.52 &           13.02 &          25.87 &    2.57 \\
\texttt{brock800\_4}  &      42.00 &       42.00 &     42.00 &         23.91 &           12.80 &          20.28 &    2.50 \\
\texttt{p\_hat300-1}  &      10.02 &       10.02 &     10.02 &         18.45 &           16.72 &           8.85 &    0.95 \\
\texttt{p\_hat300-2}  &      26.71 &       26.71 &     26.71 &        211.40 &          161.90 &          28.17 &   11.20 \\
\texttt{p\_hat300-3}  &      40.70 &       40.70 &     40.70 &         35.69 &           36.28 &          16.91 &    1.78 \\
\texttt{p\_hat500-1}  &      13.01 &       13.01 &     13.01 &         34.11 &           14.71 &          22.55 &    2.04 \\
\texttt{p\_hat500-2}  &      38.56 &       38.56 &     38.56 &        580.86 &          537.38 &          92.52 &   32.79 \\
\texttt{p\_hat500-3}  &      57.81 &       57.81 &     57.81 &         99.65 &           33.72 &          61.56 &    5.03 \\
\texttt{p\_hat700-1}  &      15.05 &       15.05 &     15.05 &         59.86 &           33.94 &          43.15 &    4.20 \\
\texttt{p\_hat700-2}  &      48.44 &       48.44 &     48.44 &       1161.67 &          295.99 &         218.26 &   71.55 \\
\texttt{p\_hat700-3}  &      71.76 &       71.76 &     71.76 &        293.90 &           93.48 &         162.79 &   15.91 \\
\texttt{p\_hat1000-1} &      17.52 &       17.52 &     17.52 &        144.76 &          119.26 &          84.75 &   10.38 \\
\texttt{p\_hat1000-2} &      54.84 &       54.84 &     54.84 &       1815.65 &          697.11 &         487.41 &  121.86 \\
\texttt{p\_hat1000-3} &      83.53 &       83.53 &     83.53 &        473.77 &          243.21 &         293.35 &   28.77 \\
\texttt{p\_hat1500-1} &      21.89 &       21.89 &     21.89 &        606.67 &          479.28 &         471.58 &   41.41 \\
\texttt{p\_hat1500-2} &      \textbf{-} &  \textbf{-} &    \textbf{76.46} &             - &               - &        1826.66 &  233.69 \\
\texttt{p\_hat1500-3} &     113.65 &      113.65 &    113.65 &       3014.42 &          879.45 &        1886.51 &  202.90 \\
\texttt{keller4}     &      13.47 &       13.47 &     13.47 &          3.35 &            1.46 &           2.69 &    0.17 \\
\texttt{keller5}     &      31.00 &       31.00 &     31.00 &        503.36 &           53.31 &         289.31 &   30.06 \\
\texttt{keller6}     &       \textbf{-} &      \textbf{63.00} &         \textbf{-} &             - &         1524.62 &               - &  146.42 \\
\texttt{sanr200\_0.7} &      23.63 &       23.63 &     23.63 &          3.44 &            1.26 &           3.21 &    0.25 \\
\texttt{sanr200\_0.9} &      48.90 &       48.90 &     48.90 &          6.04 &            1.73 &           4.80 &    0.34 \\
\texttt{sanr400\_0.5} &      20.18 &       20.18 &     20.18 &          6.60 &            3.80 &           6.19 &    0.57 \\
\texttt{sanr400\_0.7} &      33.97 &       33.97 &     33.97 &          7.21 &            4.05 &           7.49 &    0.51 \\
\texttt{MANN\_a9}     &      \textbf{17.48} &       \textbf{17.48} &     \textbf{17.47} &          0.48 &            0.28 &           0.46 &    0.05 \\
\texttt{MANN\_a27}    &     132.76 &      132.76 &    132.76 &        561.87 &            5.48 &         550.52 &   30.80 \\
\texttt{hamming6-2}  &      32.00 &       32.00 &     32.00 &          1.49 &            0.33 &           1.25 &    0.10 \\
\texttt{hamming6-4}  &       4.00 &        4.00 &      4.00 &          0.10 &            0.08 &           0.11 &    0.01 \\
\texttt{hamming8-2}  &     128.00 &      128.00 &    128.00 &        532.92 &            3.96 &         500.17 &   30.03 \\
\texttt{hamming8-4}  &      16.00 &       16.00 &     16.00 &          2.62 &            1.13 &           2.60 &    0.24 \\
\texttt{hamming10-4} &      42.67 &       42.67 &     42.67 &         97.36 &           31.77 &          93.90 &    7.46 \\
\bottomrule
\end{tabular}
}}
\caption{Results on $\vartheta_+(G)$, graphs from the second DIMACS implementation challenge.}
\label{tab:theta_ss}
\end{table}

\begin{figure}[H]
    \centering
    \includegraphics[width=0.7\textwidth]{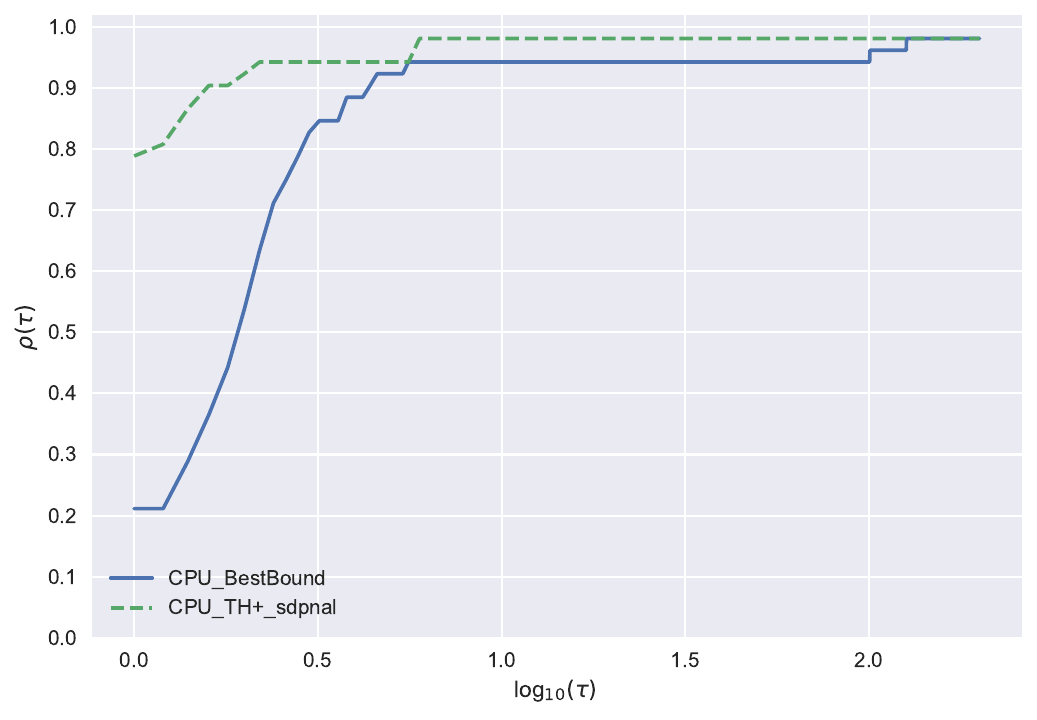}
    \caption{Performance profiles on CPU time. Comparison between BestBound and \texttt{SDPNAL+} on the computation of  $\vartheta_+(G)$.}
    \label{fig:sset}
\end{figure}

\korange{In Figure \ref{fig:sset}, we report the performance profiles obtained with respect to the CPU time needed by~\texttt{SDPNAL+} and the CPU time to identify the value of the BestBound.}
\korange{It is clear that \texttt{SDPNAL+} outperforms \texttt{ADAL-ineq} on these instances. However, we wish to highlight 
the superior} performances of~\texttt{ADAL-ineq} on the \texttt{p-hat} graphs, where we are often able to get the same bound as the optimal dual objective of~\texttt{SDPNAL+} in a much lower CPU time. 

\subsection{Comparison on instances from SDP relaxations of the coloring problem}

\korange{We now} report the results on the SDP relaxation $\bar\vartheta_+(G)$ for bounding the chromatic number of a graph. 
As before, we considered graphs from the second DIMACS implementation challenge~\cite{Johnson:1996}, available at \url{https://sites.google.com/site/graphcoloring/downloads}. 
\korange{We ran \texttt{ADAL-ineq} and \texttt{SDPNAL+}, halting the algorithms either upon satisfaction of the termination criteria with a precision of $10^{-6}$ or after a time limit of $3600$ seconds.}
\korange{As for bounding $\vartheta_+(G)$, we applied the post-processing procedure \korange{detailed in Section \ref{sec:dualbound}} every $200$ iterations and after termination of \texttt{ADAL-ineq}.}
Every dual bound computed by the post-processing procedure and, in particular, the best dual bound is a valid lower bound on the chromatic number.

\korange{In Table~\ref{tab:theta_gc1} and Table~\ref{tab:theta_gc2}, we report the same data as in Section~\ref{sec:res_clique}. As previously noted, the BestBound presented in the tables may not necessarily correspond to the result of the last post-processing call.}

\begin{figure}
    \centering
    \includegraphics[width=0.7\textwidth]{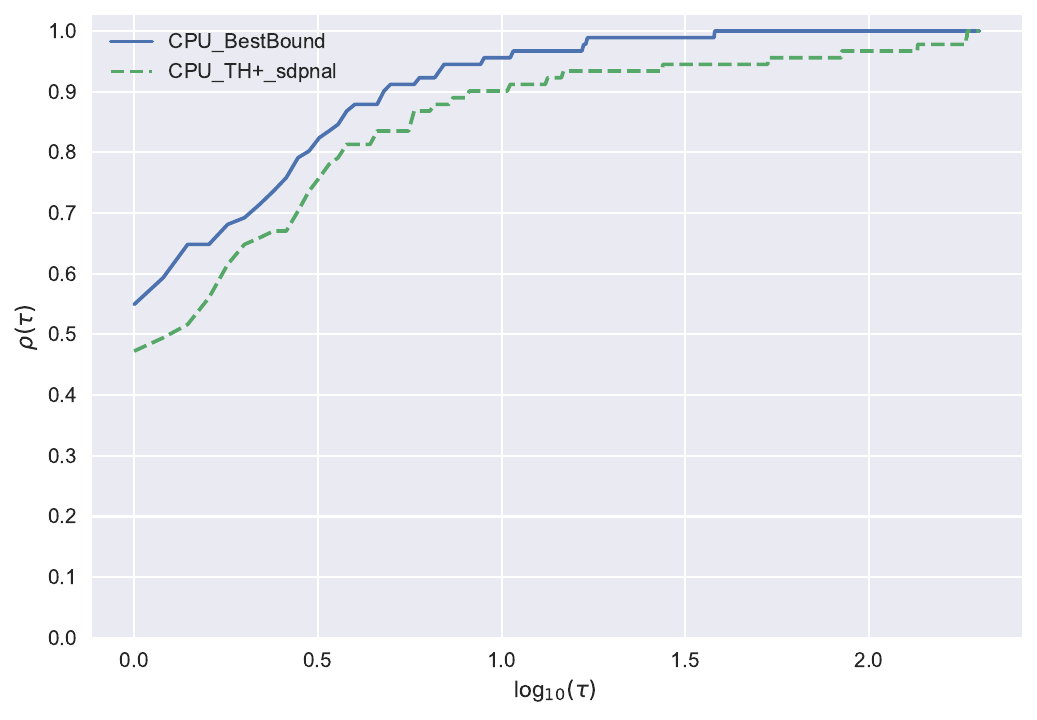}
    \caption{Performance profiles on CPU time. Comparison between BestBound and \texttt{SDPNAL+}  on the computation of $\bar\vartheta_+(G)$.}
    \label{fig:gc}
\end{figure}

We notice that the post-processing procedure fails in finding bounds on the chromatic number for several graphs \korange{for $19$ out of $113$ graphs.} This is due to the precision 
of the dual matrix given as input \korange{as it} may be too low to detect a dual feasible solution.
We also notice that, on some graphs, the bound obtained is slightly \korange{larger} than the dual objective function value obtained by~\texttt{ADAL-ineq}. \korange{This behaviour is a consequence of the precision required by Gurobi to solve the LP in the post-processing, where we set a feasibility precision of $10^{-5}$. Note that requesting a higher precision may lead to failure in the post-processing procedure. The results shown are then obtained with what, in our opinion,} is a good trade off between feasibility precision and quality of the bound.  
The CPU time needed to compute the BestBound is often much lower with respect to the time needed by~\texttt{SDPNAL+} to \korange{}{converge; this} is confirmed by the performance profiles shown in Figure~\ref{fig:gc}.
In \korange{these} profiles, we excluded \korange{the} instances for which the difference in absolute value of the BestBound found by~\texttt{ADAL-ineq} and the dual objective of~\texttt{SDPNAL+} is less than 0.5.
In particular, we excluded all the instances where the post-processing procedure was not able to compute a bound.

\begin{table}
	\centering
	{\scalebox{.65}{
	\begin{tabular}{l|rrr|rrrr}
\toprule
 &   \multicolumn{3}{c}{$\bar\vartheta_+(G)$} &  \multicolumn{4}{c}{CPU times}  \\
 \hline
Graph           &     \texttt{ADAL-ineq}       &      \texttt{SDPNAL+}       &  BestBound         &                    \texttt{ADAL-ineq}       &      \texttt{SDPNAL+}       &  BestBound             &    post-proc     \\
\midrule
\texttt{DSJC125.1}       &      4.14 &        4.14 &    4.14 &         69.01 &           22.88 &           1.72 &    0.91 \\
\texttt{DSJC125.5}       &     11.87 &       11.87 &   11.87 &          1.02 &            1.36 &           0.74 &    0.06 \\
\texttt{DSJC125.9}       &     37.80 &       37.80 &   37.80 &          1.87 &            2.03 &           1.31 &    0.09 \\
\texttt{DSJC250.1}       &      4.94 &        4.94 &    4.94 &          6.65 &            6.70 &           2.19 &    0.16 \\
\texttt{DSJC250.5}       &     16.35 &       16.35 &   16.35 &          1.90 &            2.84 &           2.00 &    0.10 \\
\texttt{DSJC250.9}       &     55.22 &       55.22 &   55.22 &          4.76 &            3.76 &           3.76 &    0.37 \\
\texttt{DSJC500.1}       &      6.25 &        6.25 &    6.25 &          8.66 &           16.57 &           5.91 &    0.16 \\
\texttt{DSJC500.5}       &     22.90 &       22.90 &   22.90 &          5.41 &            9.64 &           5.71 &    0.30 \\
\texttt{DSJC500.9}       &     84.14 &       84.14 &   84.14 &         17.04 &           16.36 &          17.53 &    1.21 \\
\texttt{DSJR500.1}       &     12.00 &       12.00 &   12.00 &         35.18 &           10.00 &          23.59 &    0.34 \\
\texttt{DSJR500.1c}      &     83.75 &       83.75 &   83.75 &             - &         1231.74 &         190.31 &  294.26 \\
\texttt{DSJR500.5}       &    \textbf{122.01} &      \textbf{122.00} &  \textbf{122.00} &        198.08 &           16.95 &         181.80 &    6.48 \\
\texttt{DSJC1000.1}      &      8.36 &        8.36 &    8.36 &         31.65 &           59.76 &          22.83 &    0.57 \\
\texttt{DSJC1000.5}      &     32.11 &       32.11 &   32.11 &         18.30 &           38.11 &          19.46 &    1.17 \\
\texttt{DSJC1000.9}      &    122.80 &      122.80 &  122.80 &         72.30 &           63.88 &          70.85 &    6.89 \\
\texttt{fpsol2.i.1}      &     65.00 &       65.00 &   65.00 &        200.57 &           11.71 &         199.60 &    2.67 \\
\texttt{fpsol2.i.2}      &     30.00 &       30.00 &   30.00 &         28.11 &            9.75 &          25.71 &    0.43 \\
\texttt{fpsol2.i.3}      &     30.00 &       30.00 &   30.00 &         27.36 &            7.82 &          27.43 &    0.42 \\
\texttt{inithx.i.1}      &     54.00 &       54.00 &   54.00 &        604.51 &           32.25 &         537.71 &    4.96 \\
\texttt{inithx.i.2}      &     \textbf{31.00} &       \textbf{31.00} &   \textbf{30.22} &        387.80 &           12.39 &          35.60 &    3.70 \\
\texttt{inithx.i.3}      &     \textbf{31.00} &       \textbf{31.00} &   \textbf{30.23} &        341.12 &           13.64 &          32.57 &    3.45 \\
\texttt{latin\_square\_10} &     \textbf{90.00} &       \textbf{89.99} &    \textbf{-} &         48.40 &           41.12 &           - &    2.91 \\
\texttt{le450\_15a}       &     \textbf{15.00} &       \textbf{15.00} &    \textbf{-} &          6.37 &            5.18 &           - &    0.12 \\
\texttt{le450\_15b}       &     15.00 &       15.00 &   15.00 &          7.06 &            5.61 &           7.13 &    0.14 \\
\texttt{le450\_15c}       &     15.00 &       15.00 &   15.00 &          3.92 &            4.70 &           4.02 &    0.10 \\
\texttt{le450\_15d}       &     15.00 &       15.00 &   15.00 &          3.86 &            4.71 &           3.96 &    0.10 \\
\texttt{le450\_25a}       &     25.00 &       25.00 &   25.00 &         19.73 &            7.54 &          19.53 &    0.29 \\
\texttt{le450\_25b}       &     \textbf{25.00} &       \textbf{25.00} &    \textbf{-} &         18.44 &            7.27 &           - &    0.23 \\
\texttt{le450\_25c}       &     25.00 &       25.00 &   25.00 &          9.67 &            7.03 &           9.78 &    0.20 \\
\texttt{le450\_25d}       &     25.00 &       25.00 &   25.00 &          9.15 &            6.82 &           9.25 &    0.19 \\
\texttt{mulsol.i.1}      &     \textbf{49.00} &       \textbf{49.00} &    \textbf{-} &         18.89 &            3.48 &           - &    0.66 \\
\texttt{mulsol.i.2}      &     31.00 &       31.00 &   31.00 &          9.02 &            2.92 &           9.04 &    0.28 \\
\texttt{mulsol.i.3}      &     31.00 &       31.00 &   31.00 &          8.13 &            3.12 &           7.47 &    0.19 \\
\texttt{mulsol.i.4}      &     31.00 &       31.00 &   31.00 &          7.97 &            2.40 &           6.31 &    0.22 \\
\texttt{mulsol.i.5}      &     31.00 &       31.00 &   31.00 &          9.88 &            3.34 &           9.01 &    0.19 \\
\texttt{school1}         &     14.00 &       14.00 &   14.00 &         14.74 &           65.03 &           8.08 &    0.40 \\
\texttt{school1\_nsh}     &     14.00 &       14.00 &   14.00 &         12.12 &           75.97 &           7.29 &    0.30 \\
\texttt{zeroin.i.1}      &     49.00 &       49.00 &   49.00 &         24.91 &            2.47 &          21.91 &    0.80 \\
\texttt{zeroin.i.2}      &     30.00 &       30.00 &   30.00 &         14.63 &            2.43 &          14.13 &    0.48 \\
\texttt{zeroin.i.3}      &     30.00 &       30.00 &   30.00 &         14.62 &            2.80 &          13.53 &    0.46 \\
\texttt{anna}            &     \textbf{11.00} &       \textbf{11.00} &    \textbf{-} &          9.97 &            1.16 &           - &    0.13 \\
\texttt{david}           &     \textbf{11.00} &       \textbf{11.00} &    \textbf{-} &          2.46 &            0.59 &           - &    0.08 \\
\texttt{huck}            &     \textbf{11.00} &       \textbf{11.00} &    \textbf{-} &          1.60 &            0.43 &           - &    0.04 \\
\texttt{jean}            &     \textbf{10.00} &       \textbf{10.00} &    \textbf{-} &          1.35 &            0.53 &           - &    0.03 \\
\texttt{games120}        &      \textbf{9.00} &        \textbf{9.00} &    \textbf{-} &          3.23 &            0.86 &           - &    0.07 \\
\texttt{miles250}        &      8.00 &        8.00 &    8.00 &          7.17 &            0.94 &           6.30 &    0.11 \\
\texttt{miles500}        &     20.00 &       20.00 &   20.00 &          6.78 &            1.69 &           6.26 &    0.11 \\
\texttt{miles750}        &     31.00 &       31.00 &   31.00 &          4.75 &            2.73 &           4.77 &    0.09 \\
\texttt{miles1000}       &     42.00 &       42.00 &   42.00 &          7.81 &            1.64 &           7.61 &    0.16 \\
\texttt{miles1500}       &     73.00 &       73.00 &   73.00 &         10.36 &            1.48 &          10.28 &    0.27 \\
\bottomrule
\end{tabular}
}}
\caption{Results on $\bar\vartheta_+(G)$, graphs from the second DIMACS implementation challenge.}
\label{tab:theta_gc1}
\end{table}

\begin{table}
	\centering
	{\scalebox{.65}{
	\begin{tabular}{l|rrr|rrrr}
\toprule
 &   \multicolumn{3}{c}{$\bar\vartheta_+(G)$} &  \multicolumn{4}{c}{CPU times}  \\
 \hline
Graph           &     \texttt{ADAL-ineq}       &      \texttt{SDPNAL+}       &  BestBound         &                    \texttt{ADAL-ineq}       &      \texttt{SDPNAL+}       &  BestBound             &    post-proc     \\
\midrule
\texttt{queen5\_5       } &      5.00 &        5.00 &    5.00 &          0.01 &            0.11 &           0.04 &    0.03 \\
\texttt{queen6\_6       } &      6.04 &        6.04 &    6.04 &          0.77 &            0.69 &           0.19 &    0.04 \\
\texttt{queen7\_7       } &      7.00 &        7.00 &    7.00 &          0.08 &            0.29 &           0.08 &    0.01 \\
\texttt{queen8\_8       } &      8.00 &        8.00 &    8.00 &          0.10 &            0.19 &           0.11 &    0.01 \\
\texttt{queen8\_12      } &     \textbf{12.00} &       \textbf{12.00} &       \textbf{-} &          0.55 &            0.62 &              - &    0.02 \\
\texttt{queen9\_9       } &      9.00 &        9.00 &    9.00 &          0.15 &            0.23 &           0.16 &    0.01 \\
\texttt{queen10\_10     } &     10.00 &       10.00 &   10.00 &          0.23 &            0.44 &           0.24 &    0.01 \\
\texttt{queen11\_11     } &     11.00 &       11.00 &   11.00 &          0.46 &            0.47 &           0.47 &    0.01 \\
\texttt{queen12\_12     } &     12.00 &       12.00 &   12.00 &          0.67 &            0.68 &           0.71 &    0.04 \\
\texttt{queen13\_13     } &     13.00 &       13.00 &   13.00 &          0.76 &            0.64 &           0.80 &    0.04 \\
\texttt{queen14\_14     } &     14.00 &       14.00 &   14.00 &          1.27 &            0.82 &           1.32 &    0.04 \\
\texttt{queen15\_15     } &     \textbf{15.00} &       \textbf{15.00} &       \textbf{-} &          1.38 &            1.26 &              - &    0.04 \\
\texttt{queen16\_16     } &     16.00 &       16.00 &   16.00 &          1.84 &            1.46 &           1.90 &    0.06 \\
\texttt{myciel3         }&      2.40 &        2.40 &    2.40 &          0.01 &            0.13 &           0.04 &    0.03 \\
\texttt{myciel4         }&      2.53 &        2.53 &    2.53 &          0.04 &            0.18 &           0.04 &    0.01 \\
\texttt{myciel5         }&      2.64 &        2.64 &    2.64 &          0.41 &            0.41 &           0.23 &    0.02 \\
\texttt{myciel6         }&      2.73 &        2.73 &    2.73 &          1.73 &            1.16 &           0.51 &    0.04 \\
\texttt{myciel7         }&      2.82 &        2.82 &    2.82 &          7.35 &            7.60 &           1.34 &    0.24 \\
\texttt{mug88\_1        } &      3.00 &        3.00 &    3.00 &         11.78 &           29.45 &           0.35 &    0.24 \\
\texttt{mug88\_25       } &      3.00 &        3.00 &    3.00 &         20.81 &           47.43 &           0.35 &    0.43 \\
\texttt{mug100\_1       } &      3.00 &        3.00 &    3.00 &         19.59 &           84.51 &           0.46 &    0.31 \\
\texttt{mug100\_25      } &      3.00 &        3.00 &    3.00 &         26.20 &           84.97 &           0.46 &    0.39 \\
\texttt{abb313GPIA      }&      8.00 &        8.00 &    8.01 &        615.04 &         2949.22 &          55.61 &    4.48 \\
\texttt{ash331GPIA      }&      3.38 &        3.38 &    3.38 &        125.34 &           17.85 &          38.24 &    1.01 \\
\texttt{ash608GPIA      }&      3.33 &        3.33 &    3.31 &        265.72 &           41.34 &         129.25 &    1.34 \\
\texttt{ash958GPIA      }&      \textbf{3.33} &        \textbf{3.33} &       \textbf{-} &        529.68 &          124.35 &           0.00 &    2.51 \\
\texttt{will199GPIA     }&      6.10 &        6.10 &    6.10 &        156.39 &           32.11 &         124.61 &    1.22 \\
\texttt{1-Insertions\_4 } &      2.23 &        2.23 &    2.23 &          1.93 &            1.01 &           0.34 &    0.07 \\
\texttt{1-Insertions\_5 } &      2.28 &        2.28 &    2.28 &         19.71 &           15.04 &           2.64 &    0.56 \\
\texttt{1-Insertions\_6 } &      2.31 &        2.31 &    2.31 &        337.22 &          100.65 &          22.80 &    3.29 \\
\texttt{2-Insertions\_3 } &      2.10 &        2.10 &    2.10 &          0.38 &            0.57 &           0.18 &    0.04 \\
\texttt{2-Insertions\_4 } &      2.13 &        2.13 &    2.13 &         25.27 &            9.06 &           1.59 &    0.36 \\
\texttt{2-Insertions\_5 } &      2.16 &        2.16 &    2.16 &        544.91 &          109.90 &          52.67 &    4.99 \\
\texttt{3-Insertions\_3 } &      2.07 &        2.07 &    2.07 &          1.22 &            1.00 &           0.31 &    0.06 \\
\texttt{3-Insertions\_4 } &      2.09 &        2.09 &    2.09 &        125.48 &           29.79 &           8.39 &    2.37 \\
\texttt{3-Insertions\_5 } &       \textbf{-}  &        \textbf{2.10} &    \textbf{2.11} &             - &         3568.47 &         130.38 &   17.94 \\
\texttt{4-Insertions\_3 } &      2.05 &        2.05 &    2.05 &          2.39 &            2.75 &           0.38 &    0.08 \\
\texttt{4-Insertions\_4 } &      2.06 &        2.06 &    2.06 &        563.58 &          130.23 &           8.89 &    6.64 \\
\texttt{1-FullIns\_3    } &      3.06 &        3.06 &    3.06 &          0.32 &            0.35 &           0.13 &    0.05 \\
\texttt{1-FullIns\_4    } &      3.12 &        3.12 &    3.12 &          4.37 &            2.45 &           1.39 &    0.10 \\
\texttt{1-FullIns\_5    } &      3.18 &        3.18 &    3.18 &         71.27 &           17.55 &          18.65 &    1.52 \\
\texttt{2-FullIns\_3    } &      4.03 &        4.03 &    4.03 &          1.34 &            0.39 &           0.74 &    0.08 \\
\texttt{2-FullIns\_4    } &      4.06 &        4.06 &    4.06 &         57.51 &            9.21 &          26.76 &    1.64 \\
\texttt{2-FullIns\_5    } &      4.08 &        4.08 &    4.08 &       2670.31 &          184.26 &         381.59 &   19.29 \\
\texttt{3-FullIns\_3    } &      5.02 &        5.02 &    5.02 &          6.12 &            1.18 &           4.47 &    0.15 \\
\texttt{3-FullIns\_4    } &      5.03 &        5.03 &    5.03 &        329.51 &           24.03 &          58.31 &    4.94 \\
\texttt{3-FullIns\_5    } &         \textbf{-} &        \textbf{5.05} &    \textbf{5.04} &             - &         1769.01 &        2965.54 &   18.40 \\
\texttt{4-FullIns\_3    } &      6.01 &        6.01 &    6.01 &         21.19 &            2.30 &           1.65 &    0.32 \\
\texttt{4-FullIns\_4    } &      6.02 &        6.02 &    6.02 &       1979.86 &           88.40 &         283.54 &   16.15 \\
\texttt{4-FullIns\_5    } &         - &           - &       - &             - &               - &              - &   14.11 \\
\texttt{5-FullIns\_3    } &      7.01 &        7.01 &    7.00 &         61.22 &            2.72 &          12.48 &    0.71 \\
\texttt{5-FullIns\_4    } &         \textbf{-} &        \textbf{7.01} &    \textbf{7.01} &       - &          207.83 &         137.49 &   21.11 \\
\texttt{wap01a          }&          \textbf{-} &       \textbf{41.00} &   \textbf{40.38} &        - &          309.86 &        3575.61 &   20.34 \\
\texttt{wap02a          }&     \textbf{40.00} &       \textbf{40.00} &    \textbf{-} &        538.62 &          473.42 &           - &    3.29 \\
\texttt{wap03a          }&     40.00 &       40.00 &   40.00 &       1594.69 &         2668.31 &        1507.80 &   10.12 \\
\texttt{wap04a          }&     40.00 &       40.00 &   40.00 &       2175.13 &         2658.54 &        2179.94 &   13.84 \\
\texttt{wap05a          }&     50.00 &       50.00 &   50.00 &       1099.11 &           24.19 &         918.93 &   11.72 \\
\texttt{wap06a          }&     \textbf{40.00} &       \textbf{40.00} &    \textbf{-} &         63.35 &           69.47 &           - &    0.74 \\
\texttt{wap07a          }&     40.00 &       40.00 &   40.00 &        309.93 &          426.97 &         145.89 &    3.00 \\
\texttt{wap08a          }&     \textbf{40.00} &       \textbf{40.00} &    \textbf{-} &        278.67 &          224.35 &                 - &    2.26 \\
\texttt{qg.order30      }&     \textbf{30.00} &       \textbf{30.00} &    \textbf{-} &         32.22 &           21.30 &                 - &    0.30 \\
\texttt{qg.order40      }&     \textbf{40.00} &       \textbf{40.00} &    \textbf{-} &        153.25 &           82.68 &                 - &    1.08 \\
\texttt{qg.order60      }&     \textbf{60.00} &       \textbf{60.00} &    \textbf{-} &       1684.60 &          496.86 &                 - &    9.74 \\
\bottomrule
\end{tabular}
}}
\caption{Results on $\bar\vartheta_+(G)$, graphs from the second DIMACS implementation challenge.}
\label{tab:theta_gc2}
\end{table}

As a further comparison between~\texttt{ADAL-ineq} and~\texttt{SDPNAL+} on SDP relaxations of the chromatic number, we built instances \korange{by} adding $1000$, $2500$ and $5000$ inequalities to $\bar\vartheta(G)$. The inequalities \korange{were} chosen randomly from those proposed by Dukanovic and Rendl~\cite{DUKANOVIC2008180} to strengthen \korange{$\bar\vartheta(G)_+$ by including to~\eqref{eq:theta_gc} the following}:
$$X_{ij} + X_{ik} - X_{jk} \leq t - 1, \quad \forall \,i, j, k \in \korange{V}. $$
\korange{Note that the main goal of this experiment is to measure the behavior of the solvers on SDPs with an increasing number of inequalities, rather than to evaluate the improvement that these valid inequalities yield over $\bar\vartheta(G)_+$.}
In Table~\ref{tab:theta_gcineq}, we report the results on some classes of graphs where 
the CPU time needed to compute the BestBound is often lower with respect to the time needed by~\texttt{SDPNAL+} to converge with a precision of $10^{-6}$.

\begin{table}
	\centering
	{\scalebox{.65}{
	\begin{tabular}{l|rrr|rrrr}
\toprule
 &   \multicolumn{3}{c}{bounds on $\chi(G)$} &  \multicolumn{4}{c}{CPU times}  \\
 \hline
Graph           &     \texttt{ADAL-ineq}       &      \texttt{SDPNAL+}       &  BestBound         &                    \texttt{ADAL-ineq}       &      \texttt{SDPNAL+}       &  BestBound             &    post-proc     \\
\midrule
\multicolumn{8}{c}{$\bar\vartheta(G)$ + $1000$ inequalities from~\cite{DUKANOVIC2008180}}\\
\midrule
DSJC500.5      &     22.90 &       22.90 &   22.90 &         15.97 &           38.73 &          13.77 &    0.54 \\
DSJC1000.1     &      8.36 &        8.36 &    8.36 &         34.43 &          154.79 &          25.85 &    0.56 \\
DSJC1000.5     &     32.11 &       32.11 &   32.11 &         33.33 &          154.82 &          34.50 &    1.18 \\
myciel7        &      2.85 &        2.85 &    2.85 &        255.84 &           13.77 &          13.39 &    6.77 \\
mug88\_25       &      \textbf{3.00} &        \textbf{3.00} &       \textbf{-} &         17.81 &           33.92 &             - &    0.25 \\
mug100\_25      &      \textbf{3.00} &        \textbf{3.00} &       \textbf{-} &         22.33 &           97.57 &             - &    0.28 \\
abb313GPIA     &      8.00 &        8.00 &    8.00 &        661.22 &         3466.57 &         178.77 &    4.49 \\
1-Insertions\_6 &      2.33 &        2.33 &    2.33 &       1001.29 &          233.15 &          65.62 &    8.98 \\
2-Insertions\_5 &      2.18 &        2.18 &    2.18 &        850.44 &          373.95 &         100.30 &    7.21 \\
3-Insertions\_5 &         \textbf{-} &           \textbf{-} &    \textbf{2.11} &             - &               - &         415.32 &   18.32 \\
4-Insertions\_4 &      2.07 &        2.07 &    2.07 &        632.60 &          363.73 &          87.32 &    6.53 \\
1-FullIns\_5    &      3.19 &        3.19 &    3.19 &       1955.97 &          141.02 &          59.14 &   35.67 \\
5-FullIns\_4    &         \textbf{-} &        \textbf{7.01} &    \textbf{7.01} &            - &          257.83 &         132.12 &   21.35 \\
wap03a         &     \textbf{40.00} &           \textbf{-} &   \textbf{40.00} &       1629.48 &               - &        1496.25 &   11.16 \\
wap04a	       &     \textbf{40.00} &	       \textbf{-} &	 \textbf{40.00} &	   2297.75 &	           - &	      2302.30 &	  13.66 \\
\midrule
\multicolumn{8}{c}{$\bar\vartheta(G)$ + $2500$ inequalities from~\cite{DUKANOVIC2008180}}  \\
\midrule
DSJC500.5      &     22.90 &       22.90 &   22.90 &         79.62 &           43.81 &          79.92 &    0.60 \\
DSJC1000.1     &      8.36 &        8.36 &    8.36 &         38.74 &          167.61 &          29.44 &    0.57 \\
DSJC1000.5     &     32.11 &       32.11 &   32.11 &        100.35 &          157.33 &          91.06 &    2.39 \\
myciel7        &      2.87 &        2.87 &    2.87 &        341.20 &           13.80 &          63.14 &    3.51 \\
mug88\_25       &      3.00 &        3.00 &    3.00 &         83.45 &           81.22 &          36.93 &    1.22 \\
mug100\_25      &      3.00 &        3.00 &    3.00 &         91.18 &           97.93 &          85.68 &    1.43 \\
abb313GPIA     &      8.00 &        8.00 &    8.00 &        685.60 &         2611.71 &         250.42 &    4.60 \\
1-Insertions\_6 &      2.34 &        2.34 &    2.34 &        827.90 &          553.82 &         111.23 &    6.91 \\
2-Insertions\_5 &      2.19 &        2.19 &    2.19 &        790.06 &          576.61 &         157.16 &    5.79 \\
3-Insertions\_5 &         \textbf{-} &        \textbf{2.13} &    \textbf{2.12} &             - &               - &        1133.88 &   18.54 \\
4-Insertions\_4 &      2.08 &        2.08 &    2.08 &        625.78 &          537.75 &         143.96 &    4.82 \\
1-FullIns\_5    &         \textbf{-} &        \textbf{3.19} &    \textbf{3.19} &             - &          145.95 &         270.64 &   34.75 \\
5-FullIns\_4    &         \textbf{-} &        \textbf{7.01} &    \textbf{7.01} &             - &          308.79 &         149.04 &   20.10 \\
wap03a	       &      \textbf{40.00} &	       \textbf{-} &	 \textbf{40.00} &	   1932.11 &	           - &	       1935.99 &   10.83 \\
wap04a	       &      \textbf{40.00} &	       \textbf{-} &	 \textbf{40.00} &	   2629.88 &	           - &	       2011.30 &   14.42 \\
\midrule
\multicolumn{8}{c}{$\bar\vartheta(G)$ + $5000$ inequalities from~\cite{DUKANOVIC2008180}}  \\
\midrule
DSJC500.5      &     22.90 &       22.90 &   22.90 &        464.82 &           46.24 &         456.20 &    1.14 \\
DSJC1000.1     &      8.36 &        8.36 &    8.36 &         63.50 &          163.05 &          48.59 &    0.70 \\
DSJC1000.5     &     32.11 &       32.11 &   32.11 &        589.69 &          168.79 &         590.91 &    2.12 \\
myciel7        &      2.89 &        2.89 &    2.89 &        583.55 &           15.44 &         132.16 &    3.90 \\
mug88\_25       &      3.00 &        3.00 &    3.00 &        199.32 &          116.73 &          78.35 &    1.56 \\
mug100\_25      &      3.00 &        3.00 &    3.00 &        216.21 &          154.84 &         190.52 &    1.39 \\
abb313GPIA     &      \textbf{8.00} &            \textbf{-} &   \textbf{8.00} &        744.60 &               - &         356.07 &    4.36 \\
1-Insertions\_6 &      2.36 &        2.36 &    2.36 &       1351.80 &          258.02 &         218.89 &    8.25 \\
2-Insertions\_5 &      2.19 &        2.19 &    2.19 &        904.85 &          572.94 &         274.81 &    5.39 \\
3-Insertions\_5 &        \textbf{-} &           \textbf{-} &    \textbf{2.12} &             - &               - &        1737.21 &   15.49 \\
4-Insertions\_4 &      2.09 &        2.09 &    2.09 &        736.69 &          831.07 &         315.23 &    4.01 \\
1-FullIns\_5    &         \textbf{-} &        \textbf{3.19} &    \textbf{3.19} &             - &          179.80 &         554.77 &   22.29 \\
5-FullIns\_4    &         \textbf{-} &        \textbf{7.01} &    \textbf{7.01} &             - &          344.28 &         293.47 &   17.54 \\
wap03a         &     \textbf{40.00} &           \textbf{-} &   \textbf{40.00} &       2591.76 &               - &        2383.76 &   15.00 \\
wap04a         &         \textbf{-} &           \textbf{-} &   \textbf{40.00} &             - &               - &        3129.61 &   18.72 \\
\bottomrule
\end{tabular}
}}
\caption{Results on bounds on $\chi(\bar G)$.}
\label{tab:theta_gcineq}
\end{table}

\section{Conclusions}\label{sec:conc}
\korange{In this paper, we} propose a numerical comparison between \korange{\texttt{ADAL-ineq}, an enhanced version of \texttt{ADAL} where the presence of linear inequality constraints is smartly handled,} and \texttt{SDPNAL+}, the state-of-the-art method for solving large-scale SDPs \korange{that has been awarded} the Beale-Orchard-Hays Prize in 2018.
We consider random instances as well as instances from the SDP relaxations of the graph coloring 
problem and the maximum clique problem. 
The post-processing procedure used is developed to obtain a dual feasible solution
which in turn gives a bound on the optimal primal value.
From a practical \korange{standpoint}, as long as we use SDPs to address combinatorial optimization problems, the post-processing procedure allows to stop the execution of the ADMM
as soon as a ``good" bound is obtained, even if the convergence criterion is far \korange{from being} met.
Furthermore, the fact that a dual feasible solution is detected, allows to use re-optimization techniques
within branch-and-bound frameworks \korange{and is} what we plan to focus on in the near future.

\section*{Acknowledgements}
The authors acknowledge support within the project RM120172A2970290
which has received funding from Sapienza, University of Rome.
They are also indebted to Fabrizio Rossi and Stefano Smriglio for their useful suggestions and to Griffin D. Kent for his precious help in proof reading this manuscript.
Last but not least, they are grateful to two anonymous referees for the careful reading of the manuscript and the valuable comments that helped to improve the paper.

\section*{Declarations}
\paragraph{Conflict of interest} The authors declare no conflict of interest

\clearpage
\bibliographystyle{plain}
\bibliography{main}

\end{document}